\documentclass{amsart}
\usepackage{graphicx}
\usepackage{amssymb}
\usepackage{epstopdf}
\usepackage{nicefrac}
\usepackage{enumerate}
\usepackage{float}
\usepackage{color}
\usepackage{pst-all}
\usepackage{tabularx}
\usepackage{hyperref}
\usepackage{caption}
\usepackage{subcaption}

\newtheorem{thm}{THEOREM}[section]

\newtheorem{cor}[thm]{COROLLARY}
\newtheorem{defn}[thm]{DEFINITION}

\newtheorem{hyp}[thm]{HYPOTHESIS}
\newtheorem{lemma}[thm]{LEMMA}
\newtheorem{prob}[thm]{PROBLEM}
\newtheorem{prop}[thm]{PROPOSITION}

\newtheorem*{prob*}{PROBLEM}

\newcommand{\ds}{\displaystyle}

\newcommand{\e}{{\epsilon}}
\newcommand{\g}{{\gamma}}
\newcommand{\wmW}{\widehat{\mW}}
\newcommand{\bRt}{{\bf R}_{0}}

\newcommand{\bR}{{\bf R}}

\newcommand{\mD}{{\mathbb D}}
\newcommand{\mK}{{\mathbb K}}
\newcommand{\mN}{{\mathbb N}}
\newcommand{\mR}{{\mathbb R}}
\newcommand{\mS}{{\mathbb S}}
\newcommand{\mT}{{\mathbb T}}
\newcommand{\mZ}{{\mathbb Z}}
\newcommand{\mW}{{\mathbb W}}
\newcommand{\cA}{{\mathcal A}}
\newcommand{\cC}{{\mathcal C}}
\newcommand{\cD}{{\mathcal D}}

\newcommand{\cG}{{\mathcal G}}
\newcommand{\cI}{{\mathcal I}}
\newcommand{\cK}{{\mathcal K}}
\newcommand{\cL}{{\mathcal L}}

\newcommand{\cO}{{\mathcal O}}

\newcommand{\cR}{{\mathcal R}}

\newcommand{\cW}{{\mathcal W}}
\newcommand{\cX}{{\mathcal X}}
\newcommand{\cZ}{{\mathcal Z}}

\newcommand{\cGM}{\cG_{\fM}}

\newcommand{\fC}{{\mathfrak{C}}}
\newcommand{\fD}{{\mathfrak{D}}}
\newcommand{\fM}{{\mathfrak{M}}}

\newcommand{\fU}{{\mathfrak{U}}}

\newcommand{\fZ}{{\mathfrak{Z}}}
\newcommand{\fW}{{\mathfrak{W}}}

\newcommand{\oU}{\overline{U}}

\newcommand{\otheta}{\overline{\theta}}
\newcommand{\ovg}{\overline{\gamma}}
\newcommand{\ovl}{\overline{\lambda}}

 \newcommand{\vp}{{\varphi}}

\newcommand{\A}{{\rm Area}}

\begin{document}

\title{Perspectives on Kuperberg flows}

  \dedicatory{Dedicated to Professor Krystyna Kuperberg  on the occasion
     of her $70$th birthday}

\author{Steven Hurder}
\address{Steven Hurder, Department of Mathematics, University of Illinois at Chicago, 322 SEO (m/c 249), 851 S. Morgan Street, Chicago, IL 60607-7045}
\email{hurder@uic.edu}
\thanks{Preprint date: July 4, 2016}

\author{Ana Rechtman}
\address{Ana Rechtman, Institut de Recherche Math\'ematique Avanc\'ee,
Universit\'e de Strasbourg,
7 rue Ren\'e Descartes,
67084 Strasbourg, France}
\email{rechtman@math.unistra.fr}

\date{}


% \keywords{}
 
 \begin{abstract}
The ``Seifert Conjecture'' asks, ``Does   every non-singular vector field on the 3-sphere $\mS^3$ have a periodic orbit?''  
In a celebrated work, Krystyna Kuperberg   gave a construction of a smooth aperiodic vector field on a plug, which is then used to construct    counter-examples to the  Seifert Conjecture for smooth flows on the $3$-sphere, and on compact 3-manifolds in general. The dynamics of the flows in these plugs have been extensively studied, with more precise results known  in special ``generic'' cases of the construction. Moreover, the dynamical properties of smooth perturbations of Kuperberg's construction have been considered.  
    In this work, we recall some of the results obtained to date for the Kuperberg flows and their perturbations. Then the main point of this work is to focus attention on how the known results for Kuperberg flows depend on the assumptions imposed on the flows, and to discuss some of the many interesting questions and problems that remain open about their dynamical and ergodic properties.
 \end{abstract}

\maketitle

  \tableofcontents

 \vfill
 \eject

\section{Introduction} \label{sec-intro}

The ``Seifert Conjecture'', as originally formulated in 1950 by Seifert   \cite{Seifert1950},  asked: ``Does   every non-singular vector field on the $3$-sphere $\mS^3$ have a periodic orbit?''  
This problem became more specific  following the construction in 1966 by  F.W.~Wilson     \cite{Wilson1966}    of a smooth flow on a  \emph{plug}  with exactly two periodic orbits. (See also the later paper  with Percel     \cite{PercelWilson1977}.) These Wilson plugs  could then be used to modify any given flow on a $3$-manifold to obtain one with only isolated periodic orbits. Thus,   the Seifert Conjecture is reduced to showing that a flow on a $3$-manifold with only a finite number of periodic orbits can be perturbed to one with no periodic orbits.

Paul Schweitzer   showed in the  1974 work \cite{Schweitzer1974}   how to modify the construction of the Wilson plug in several fundamental ways, to obtain a $C^1$-flow with no periodic orbits in the plug.  (See also Rosenberg's survey   \cite{Rosenberg1974}.)   He then used this plug to show that for any closed   $3$-manifold $M$, there exists  a non-singular $C^1$-vector field on $M$ without periodic orbits.  Jenny Harrison  constructed in  \cite{Harrison1988} a   modified  version of the Schweitzer construction,  which she used to construct  an aperiodic  plug with   a $C^2$-flow. Finally, in 1994 
Krystyna Kuperberg  constructed in \cite{Kuperberg1994}   counter-examples   to the smooth version of the Seifert Conjecture.
 \begin{thm}  [Kuperberg] \label{thm-mainK}
On every closed oriented 3-manifold $M$, there exists  a $C^\infty$  non-vanishing vector field without periodic orbits.
\end{thm}
The fascination with this result has many sources. The proof itself introduced  a construction  of aperiodic\linebreak $3$-dimensional plugs which was   notable for its simplicity and beauty.  
The   Kuperberg construction  remains the only general method to date
to construct  $C^\infty$-flows without periodic orbits on arbitrary
closed 3-manifolds.  The introductory comments in the S{\'e}minaire
Bourbaki article by Ghys \cite{Ghys1995} discusses the background  and obstacles to the construction of aperiodic flows. Kuperberg has given a brief introduction to the use of plugs in her proof of this result here 
\cite{Kuperberg1999}.

The dynamical properties of the Kuperberg flows in these plugs also
have many special qualities. The flows must have   zero topological
entropy, which follows by an observation of Ghys in  \cite{Ghys1995}
that an aperiodic flow on a 3-manifold has topological entropy equal to zero, that is
a consequence of a well-known result of Katok in \cite{Katok1980}.  The aperiodic flows in the Kuperberg plugs have a unique minimal set, whose structure is unknown in general, and moreover the generic Kuperberg  flow   preserves a $2$-dimensional compact lamination with boundary that is contained in the interior of the plug, and contains the minimal set. Finally, 
the Kuperberg flows have an open set of non-recurrent points that limit on the minimal set, by a result of Matsumoto \cite{Matsumoto1995}, which adds more to the complexity of the dynamical properties of these flows.

   There followed after Kuperberg's seminal work,  a collection of    works explaining in further detail the proof of the aperiodicity for the   Kuperberg flow in a plug, and investigated its dynamical properties: 
  \begin{itemize}
\item the S{\'e}minaire Bourbaki lecture  \cite{Ghys1995} by \'{E}tienne Ghys; 
\item  the notes  by  Shigenori Matsumoto \cite{Matsumoto1995} in Japanese,     later translated into English; 
\item  the joint paper     \cite{Kuperbergs1996} by Greg Kuperberg and Krystyna Kuperberg; 
\item  the monograph     \cite{HR2016a} by the authors. 
\end{itemize} 
 Moreover, a notable feature about the construction of these plugs, is that there are many choices in their construction, which  influence the global dynamics of the resulting flows. As Ghys wrote in \cite[page 302]{Ghys1995}: 
\begin{quote}
\emph{Par ailleurs, on peut construire beaucoup de pi\`eges de Kuperberg et il n'est pas clair qu'ils aient le m\^eme dynamique.}
\end{quote}

The works cited above suggest many interesting open problems
concerning the dynamics of Kuperberg flows. Also,    the dynamical properties of flows constructed by variants of the Kuperberg construction which are not aperiodic are studied in  \cite{HR2016b}, showing that the Kuperberg flows are indeed very special, as they lie at the ``boundary of chaos'' in the $C^{\infty}$-topology on flows.

 The purpose of this paper is to collect together these open problems, as well as to formulate more precisely further questions about Kuperberg flows which have arisen in their study, to illuminate the surprising   dynamical complexities of the flows that can arise in this  very special class of examples.

We first  describe the construction of the Kuperberg plugs in Sections~\ref{sec-wilson} and \ref{sec-kuperberg}, with an emphasis on the choices involved.    Section~\ref{sec-wilson} gives the construction of the modified Wilson plugs, which are the foundation of the construction.    Section~\ref{sec-kuperberg}   gives the   construction of the Kuperberg plugs. These constructions are given in a succinct manner, and the interested reader should consult the literature cited above for further details and discussions.

   Section~\ref{sec-generic}  introduces the additional assumptions imposed on the constructions of Kuperberg flows in the works \cite{HR2016a,HR2016b}, which are called the \emph{generic hypotheses}.  We also  introduce variations of these generic hypotheses, whose implications for the dynamics of the flows    will be discussed in later sections.

 It was observed in the  works \cite{Ghys1995,Matsumoto1995}   that
 any orbit not escaping the plug in forward or backward time, limits
 to the invariant set defined as the closure of the ``special orbits''
 for the flow. One consequence is that the Kuperberg flow in a plug has   a unique
 minimal  set, denoted by
 $\Sigma$. Theorem~\ref{thm-minimal} recalls these results. It is a
 remarkable aspect of the construction of the Kuperberg flows, that
 they preserve an embedded infinite surface with boundary, denoted by $\fM_0$, which contains the special orbits, and so that its closure $\fM = \overline{\fM_0}$ is a type of lamination with boundary that contains the unique minimal set. The relation between the two sets $\Sigma \subset \fM$ is an important  theme in the study of the dynamics of Kuperberg flows.  Section~\ref{sec-lamination} gives an outline of the structure theory for   the embedded surface $\fM_0$ in the case of generic flows, which was developed in \cite{HR2016a}. This structure theory, and the corresponding properties of the level function defined on $\fM_0$, are crucial for the proofs of many of the properties of the minimal set $\Sigma$ in \cite{HR2016a}.

 In  Section~\ref{sec-denjoy}, we consider the relation between  the
 minimal set $\Sigma$ and the laminated space $\fM$, and recall the
 conditions necessary to show they are equal. This result is a type of
 ``Denjoy Theorem'' for laminations,  and its proof in \cite{HR2016a}
 relies on the generic hypotheses on the flow in fundamental ways. 
 We consider  how the properties of the embedding
 $\Sigma \subset \fM$ depend on variations of the  generic hypotheses
 imposed on the flow. 

An interesting problem is a general formulation of a Denjoy Theorem
 for laminations, that we present here as it is independent of the
 construction of the Kuperberg plug.

\begin{prob*}[Problem \ref{prob-denjoy}]
Let $\cL$ be a compact connected 2-dimensional and codimension 1
lamination, and let $\cX$ be a smooth vector field tangent to the
leaves of $\cL$. If $\cL$ is minimal and the flow of $\cX$ has no
periodic orbits, show that every orbit is dense.
\end{prob*}
 
Ghys observed in \cite{Ghys1995} that an aperiodic flow must have entropy equal to zero, using a well-known result of Katok \cite{Katok1980}, and thus the Kuperberg flows must have zero entropy.  In the work \cite{HR2016a}, the authors showed that in fact, while the usual entropy of the flow vanishes, the ``slow entropy'' as defined in \cite{deCarvalho1997,KatokThouvenot1997} of a generic Kuperberg flow is positive for exponent $\alpha = 1/2$. This calculation used the fact that the embedded surface $\fM_0 \subset \mK$ has subexponential but not polynomial growth rate, which follows from the structure theory developed for it in the generic case. Finally, these results suggests the study of the Hausdorff dimensions of the sets $\Sigma$ and $\fM$ and how they depend on the choices used in the construction of  the Kuperberg flow. These and related questions are addressed in  Section~\ref{sec-entropy}.

A standard problem in topological dynamical systems theory is to
describe the topological type of the closed attractors for the system, and for
closed invariant transitive subsets more generally. Attractors often
have very complicated topological description, and the theory of shape
for spaces \cite{MardesicSegal1982,Mardesic1999} is used to describe
them. For example, the shape of the unique minimal set $\Sigma$ for a
\emph{generic} Kuperberg flow is shown in
\cite{HR2016a} to be \emph{not stable}, but to satisfy a \emph{Mittag-Leffler Property} on its
homology groups. The proof of these assertions requires essentially
all of the material developed in the monograph
\cite{HR2016a}. Describing  the shape of a  dynamically defined
invariant set of an arbitrary  flow is typically quite difficult, but
also can be highly revealing about the dynamical properties of the
flow.  In Section~\ref{sec-shape}, we discuss questions concerning
the shape properties of the minimal set for general Kuperberg flows.

The proofs in \cite{HR2016a} suggest a strong relation between the
shape approximations of the minimal set $\Sigma$ and the entropy of
the flow. The following problem can be stated for general flows and
the motivation is given in Section~\ref{sec-shape} and Problem~\ref{prob-shape4}.

\begin{prob*}
Assume that a flow has an exceptional minimal set whose shape is
 not stable. Is the slow entropy of the flow positive?
\end{prob*}

A minimal set is said to be exceptional if it is not a submanifold of
the ambient manifold.
We conjecture that part of the hypothesis needed to show that the
shape is related with the entropy of the flow is that the minimal set
has ``small'' dimension, that is dimension 1 or 2. 

The Derived from Kuperberg flows, or DK--flows, were introduced in
\cite{HR2016b}, and are obtained by varying the construction of the
usual Kuperberg flows so that the periodic orbits are not ``broken
open''. Thus, the DK--flows are quite useless as counterexamples to
the Seifert Conjecture, but they are obtained by smooth variations of
the standard Kuperberg flows, so are of central  interest from the point of view
of the properties of Kuperberg flows in the space of flows. The work
\cite{HR2016b} gave constructions of DK--flows which in fact have
countably many independent horseshoe subsystems, and thus have
positive topological entropy.  Moreover, these examples can be
constructed arbitrarily close to the   generic Kuperberg flows.
Section~\ref{sec-DK} discusses a variety of questions about the flows
which are $C^{\infty}$-close to Kuperberg flows. One topic in
particular is notable, that the horseshoes  generated by a variation
of the construction, are created using the shape approximations
discussed in Section~\ref{sec-shape}, providing more reasons to
explore the relation between shape and entropy for flows.

The authors dedicate this work to Krystyna Kuperberg, both for her discovery of the class of dynamical systems   introduced in her celebrated works on aperiodic flows, and whose 
  comments and suggestions to the authors have been important both  over the course of   writing   the monograph \cite{HR2016a}, and have inspired our continued fascination with  ``Kuperberg flows''.

\section{Modified Wilson plugs}\label{sec-wilson}

In this section and the next, we present the construction of the  Kuperberg plugs which are the basis for the proof of Theorem~\ref{thm-mainK}, with commentary on the choices made in the process.   First, we recall that a ``plug'' is a manifold with boundary endowed with a flow, that
 enables  the modification of a given flow on a $3$-manifold inside a flow-box. The idea is that 
 after modification by insertion of a plug, a periodic orbit for the given flow is ``broken
 open'' -- it   enters the plug and never
 exits. Moreover,   Kuperberg's construction  does this modification without introducing additional  periodic orbits.
 The first step is to construct Kuperberg's modified Wilson plug, which is analogous to the modified Wilson plug used by Schweitzer in     \cite{Schweitzer1974}.

The notion of a ``plug'' to be inserted in a flow on a $3$-manifold was introduced by Wilson \cite{Wilson1966,PercelWilson1977}.  
A $3$-dimensional plug is a manifold $P$ endowed with a vector field $\cX$ satisfying the following conditions. The 3-manifold $P$ is of the form $D \times [-2,2]$, where $D$ is a compact 2-manifold with boundary $\partial D$. Set 
$$\partial_v P = \partial D \times [-2,2] \quad , \quad \partial_h^- P = D \times \{-2\} \quad , \quad \partial_h^+ P = D \times \{2\} \ .$$
Then   the boundary  of $P$ has a decomposition
$$\partial P ~ = ~  \partial_v P \cup \partial_h P ~ = ~  \partial_v P \cup \partial_h^-P \cup \partial_h^+ P \ .$$
Let $\frac{\partial}{\partial z}$ be the \emph{vertical} vector field on $P$, where $z$ is the coordinate of the interval $[-2,2]$.

The vector field $\cX$ must satisfy the conditions:
\begin{itemize}
\item[(P1)] \emph{vertical at the boundary}: $\cX=\frac{\partial}{\partial z}$ in a neighborhood of $\partial P$; thus, $\partial_h^- P$ and $\partial_h^+ P$  are the entry and exit regions of $P$ for the flow of $\cX$, respectively; 
\item[(P2)]   \emph{entry-exit condition}: if a point $(x,-2)$ is in the same trajectory as $(y,2)$, then $x=y$. That is,  an orbit that traverses $P$, exits just in front of its entry point;
\item[(P3)] \emph{trapped orbit}: there is at least one entry point whose entire \emph{forward} orbit is contained in $P$; we will say that its orbit is \emph{trapped} by $P$;
\item[(P4)] \emph{tame}: there is an embedding $i \colon P\to \mR^3$ that preserves the vertical direction.
\end{itemize}

   Note that conditions (P2) and (P3) imply that if the forward orbit of a point $(x,-2)$ is trapped, then the backward orbit of $(x,2)$ is also trapped. 

A {\it semi-plug} is  a manifold $P$ endowed with a vector field $\cX$ as above, satisfying conditions (P1), (P3) and (P4), but not necessarily (P2).  The concatenation of a semi-plug with an inverted copy of it, that is a copy where the direction of the flow is inverted, is then a plug. 

Note that condition   (P4)  implies that given any open ball  $B(\vec{x},\e) \subset \mR^3$ with $\e > 0$, there exists a modified embedding $i' \colon P \to B(\vec{x},\e)$ which preserves the vertical direction again.   Thus, a plug can be used to change a vector field $\cZ$ on any $3$-manifold $M$ inside a flowbox, as follows. Let $\vp \colon U_x \to (-1,1)^3$ be a coordinate chart which maps the vector field $\cZ$ on $M$ to the vertical vector field $\frac{\partial}{\partial z}$. Choose a modified embedding $i' \colon P \to B(\vec{x},\e) \subset (-1,1)^3$, and then  replace the flow $\frac{\partial}{\partial z}$ in the interior of $i'(P)$ with the image of $\cX$. This results in a   flow $\cZ'$ on $M$.

 The   entry-exit condition implies that    a  periodic orbit of $\cZ$ which meets $\partial_h P$ in a non-trapped point,   will remain periodic after this modification. An orbit of $\cZ$ which meets $\partial_h P$ in a  trapped point  never exits the plug $P$, hence  after modification,   limits to a closed invariant   set contained in $P$.  A closed invariant set contains a minimal set for the flow, and thus, a plug  serves as a device    to insert a minimal set  into a flow. 

In the work of Wilson \cite{Wilson1966}, the basic plug has the shape
of a solid cylinder, whose base $\ds \partial_h^- P = D \times \{-2\}$
is a planar disk $\mD^2$. Schweitzer introduced in \cite{Schweitzer1974} plugs for which the base is obtained from a $2$-torus minus an open disk, so $\ds \partial_h^- P \cong  {\mT^2}   - \mD^2$, which  has the homotopy type of a wedge of two circles. As we shall see below, a key idea behind the Kuperberg construction is to consider a base which is obtained from an annulus by adding two connecting strips, so has the homotopy type of three circles.
 The     ``modified Wilson Plug''   is a flow on a cylinder minus its core.  The flow   has  two periodic orbits, and the dynamics of the flow is not   stable under perturbations. The  instability of its dynamics is a key property for  its role in breaking open periodic orbits.

   The first step in the construction is to define a flow on a rectangle as follows.  
The rectangle is defined by
\begin{equation}\label{eq-rectangle}
{\bR} = [1,3]\times[-2,2] = \{(r,z) \mid 1 \leq r \leq 3 ~ \& -2 \leq z \leq 2\} \ .
\end{equation}
For a constant $0 < g_0\leq 1$,  choose a $C^\infty$-function $g \colon \bR \to [0,g_0]$   which satisfies the ``vertical'' symmetry condition $g(r,z) = g(r,-z)$. 
Also, require that   $g(2,-1) = g(2,1) = 0$,   that $g(r,z) = g_0$ for $(r,z)$ near the boundary of $\bR$, and that $g(r,z) > 0$ otherwise. 

\begin{figure}[!htbp]
\centering
 \includegraphics[width=30mm]{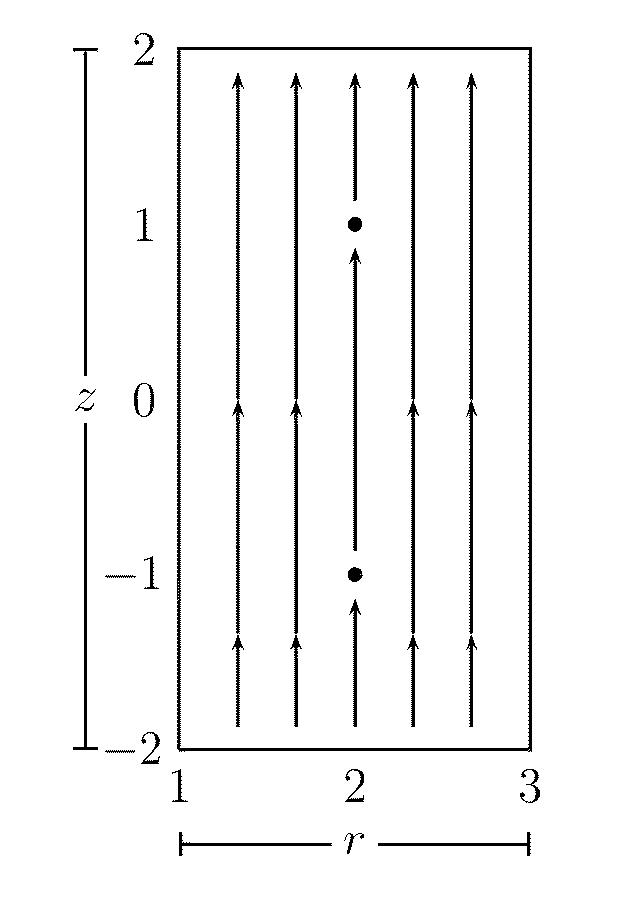}
\caption{Vector field $\cW_v$ \label{fig:wilson1}}
\vspace{-6pt}
\end{figure}

Define the vector field $\cW_v = g \cdot \frac{\partial}{\partial  z}$ which has two singularities, $(2,\pm 1)$, and   is otherwise everywhere vertical. The flow lines of this vector field are    illustrated in Figure~\ref{fig:wilson1}.  The value of $g_0$ chosen influences the quantitative nature of the flow, as small values of $g_0$ result in a slower vertical climb for the flow, but does not alter the qualitative nature of the flow.    

The next step is to suspend  the flow of the vector field $\cW_v$ to obtain a flow on the $3$-dimensional plug. This is done as follows, where we make more precise choices of the suspension flow than given in  \cite{Kuperberg1994}, though these choices do not matter so much. Choose a $C^\infty$-function $f \colon \bR \to [-1,1]$ which satisfies the following conditions:
\begin{enumerate}
\item[(W1)] $f(r,-z) = -f(r, z)$ ~ [\emph{anti-symmetry in z}]
\item[(W2)] $f(r,z) = 0$ for $(r,z)$ near the boundary of $\bR$
\item[(W3)] $f(r,z) \geq 0$ for  $-2 \leq z \leq 0$. 
\item[(W4)] $f(r,z) \leq 0$ for  $0  \leq z \leq 2$. 
\item[(W5)] $f(r,z) =1$ for $5/4 \leq r \leq 11/4$  ~ \text{and} ~  $-7/4 \leq z \leq -1/4$. 
\item[(W6)] $f(r,z) = -1$ for $5/4 \leq r \leq 11/4$ ~ \text{and}~  $1/4 \leq z \leq 7/4$. 
\end{enumerate}
Condition (W1) implies that $f(r,0) =0$ for all $1 \leq r \leq 3$. The other conditions (W2), (W3), and (W4) are assumed in the works 
 \cite{Ghys1995,Kuperberg1994,Matsumoto1995} while (W5) and (W6) were imposed in  \cite{HR2016a} in order to facilitate the description of the dynamics of the Kuperberg flows, but do not qualitatively change the resulting dynamics.

Now define the manifold with boundary
\begin{equation}\label{eq-wilsoncylinder}
\mW=[1,3] \times \mS^1\times[-2,2] \cong {\mathbf R} \times \mS^1
\end{equation}
with cylindrical coordinates  $x = (r, \theta,z)$.  That is, 
 $\mW$ is a      solid cylinder with an open core removed, obtained by rotating  the rectangle $\bR$,  considered as embedded in $\mR^3$, around the $z$-axis. 
 
Extend the functions $f$ and $g$ above to $\mW$ by setting $f(r, \theta, z) = f(r, z)$ and $g(r, \theta, z) = g(r, z)$, so that they are   invariant under rotations around the $z$-axis. 
The modified  Wilson vector field on $\mW$  is given  by 
\begin{equation}\label{eq-wilsonvector}
\cW =g(r, \theta, z)  \frac{\partial}{\partial  z} + f(r, \theta, z)  \frac{\partial}{\partial  \theta} \ .
\end{equation}
Observe that the vector field $\cW$ is vertical near the boundary of $\mW$ and horizontal in the periodic orbits. 
Also, $\cW$ is tangent to the cylinders $\{r=cst\}$. 

\begin{figure}[!htbp]
\centering
\begin{subfigure}[c]{0.3\textwidth}{\includegraphics[height=28mm]{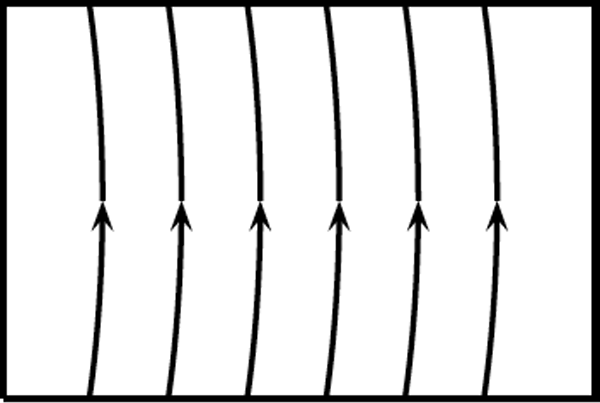}}\caption{$r \approx 1,3$}\end{subfigure}
\begin{subfigure}[c]{0.3\textwidth}{\includegraphics[height=28mm]{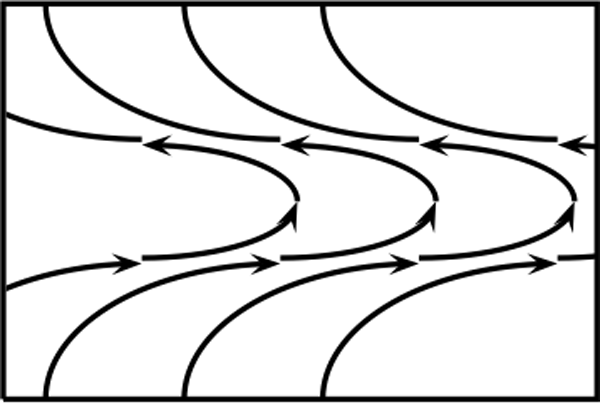}}\caption{$r \approx 2$}\end{subfigure}
\begin{subfigure}[c]{0.3\textwidth}{\includegraphics[height=28mm]{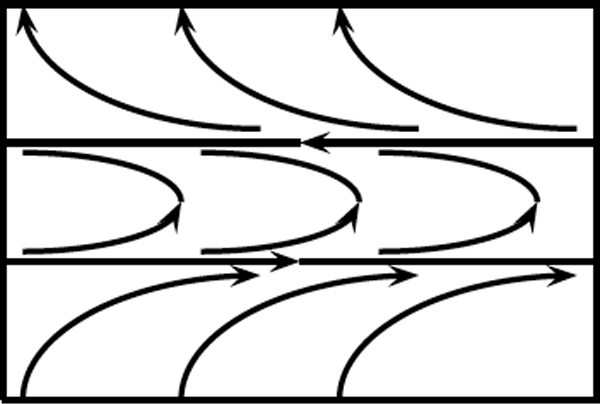}}\caption{$r =  2$}\end{subfigure}
\caption{$\cW$-orbits on the cylinders $\{r=const.\}$ \label{fig:flujocilin}}
\vspace{-6pt}
\end{figure}

\begin{figure}[!htbp]
\centering
\begin{subfigure}[c]{0.4\textwidth}{\includegraphics[height=54mm]{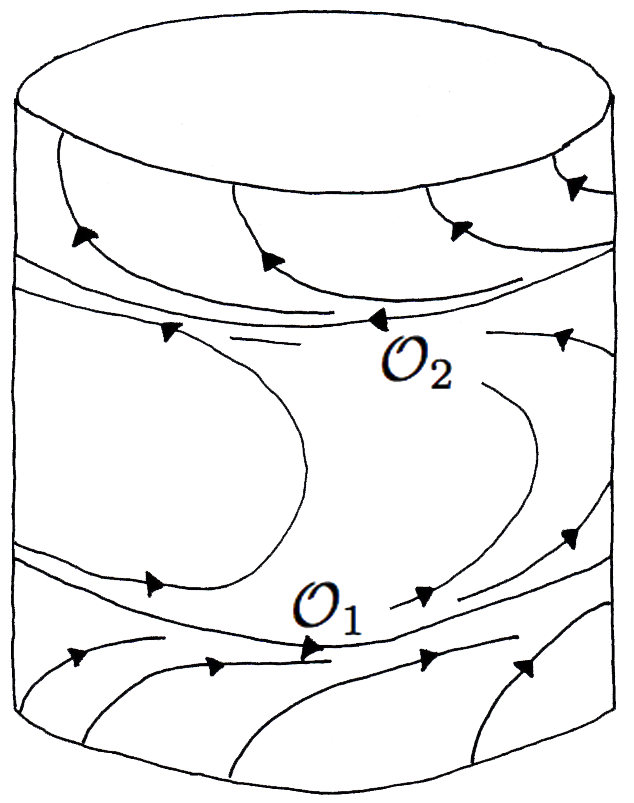}}\end{subfigure}
\begin{subfigure}[c]{0.4\textwidth}{\includegraphics[height=54mm]{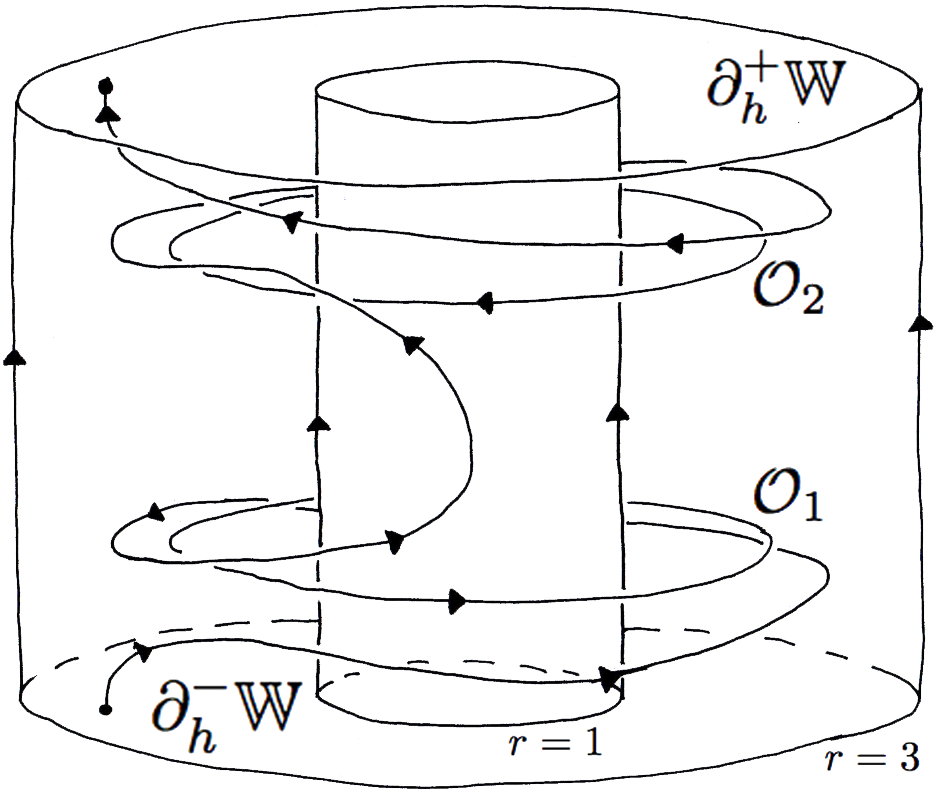}}\end{subfigure}
 \caption[{$\cW$-orbits in the cylinder  $\cC=\{r=2\}$}]{$\cW$-orbits in the cylinder  $\cC=\{r=2\}$ and in $\mW$ \label{fig:Reebcyl}}  
%\vspace{-6pt}
\end{figure}

Let $\Psi_t$ denote the flow of $\cW$ on $\mW$.  
The flow  of $\Psi_t$ restricted to the cylinders $\{r=cst\}$ is illustrated (in cylindrical coordinate slices) by  the lines in  Figures~\ref{fig:flujocilin} and \ref{fig:Reebcyl}.   The flow  of $\Psi_t$ restricted to the cylinder $\{r=2\}$ in Figure~\ref{fig:flujocilin}(C) is a called the \emph{Reeb flow}, which Schweitzer remarks in \cite{Schweitzer1974} was the inspiration for his introduction of this variation on the Wilson plug.

  We will    make   reference to the following sets in $\mW$:
  \begin{eqnarray*}
\cC ~ & \equiv & ~  \{r=2\}   \quad \text{[\emph{The ~ Full ~ Cylinder}]}\\
\cR ~ & \equiv & ~  \{(2,\theta,z) \mid  -1 \leq z \leq 1\} \quad \text{[\emph{The ~ Reeb ~ Cylinder}]}\\
\cA ~ & \equiv & ~  \{z=0\} \quad \text{[\emph{The ~ Center ~ Annulus}]}\\
\cO_i ~ & \equiv & ~  \{(2,\theta,(-1)^i)  \} \quad \text{[\emph{Periodic Orbits, i=1,2}]}
\end{eqnarray*}
Note that  $\cO_1$ is the lower   boundary circle  of the Reeb cylinder $\cR$, and $\cO_2$ is the upper boundary circle.

 Let us also recall some of the  basic properties of the modified  Wilson flow, which follow from the construction of the vector field $\cW$ and the conditions (W1) to (W4)    on the suspension function  $f$.
 
Let $R_{\varphi} \colon \mW \to \mW$ be rotation by the angle $\varphi$. That is, $R_{\varphi}(r,\theta,z) = (r, \theta + \varphi, z)$.  
\begin{prop}\label{prop-wilsonproperties}
Let $\Psi_t$ be the flow on $\mW$ defined above, then:
\begin{enumerate}
\item $R_{\varphi} \circ \Psi_t = \Psi_t \circ R_{\varphi}$ for all $\varphi$ and $t$.
\item The flow $\Psi_t$ preserves the cylinders $\{r=const.\}$   and  in particular 
    preserves the cylinders $\cR$ and $\cC$. 
\item $\cO_i$ for $i=1,2$ are  the periodic orbits for $\Psi_t$.
\item For $x = (2,\theta,-2)$, the forward orbit   $\Psi_t(x)$ for $t > 0$ is trapped.
\item For $x = (2,\theta,2)$, the backward orbit  $\Psi_t(x)$ for $t < 0$ is trapped.
\item For $x = (r,\theta,z)$ with $r \ne 2$, the orbit $\Psi_t(x)$ terminates in the top face $\partial_h^+ \mW$ for some $t \geq 0$, and terminates in $\partial_h^- \mW$ for some $t \leq 0$.
\item The flow $\Psi_t$  satisfies   the entry-exit condition {\rm (P2)} for plugs.
\end{enumerate}
\end{prop}
 The properties of the flow $\Psi_t$ on $\mW$ given in Proposition~\ref{prop-wilsonproperties} are fundamental for showing that the Kuperberg flows constructed in the next section are aperiodic. On the other hand, the study of the further dynamical properties of the Kuperberg flows reveals the importance of the behavior of the flow $\Psi_t$ in open neighborhoods of the periodic orbits $\cO_1$ and $\cO_2$. This behavior depends strongly on the properties of the function $g$ in open neighborhoods of its vanishing points $(2, \pm 1)$, as will be discussed in later sections.  In particular, note that if the function $g$ is modified in arbitrarily small neighborhoods of the points $(2,-1)$ and $ (2,1)$,  so that $g(r,z) > 0$ on $\bR$, then the resulting flow on $\mW$ will have no periodic orbits, and no trapped orbits.

\section{The Kuperberg plugs}\label{sec-kuperberg}
 
Kuperberg's construction in \cite{Kuperberg1994} of   aperiodic  smooth flows on plugs   introduced a fundamental new idea,   that of  ``geometric surgery'' on the modified     Wilson plug $\mW$ constructed in the previous section, to obtain  the \emph{Kuperberg Plug} $\mK$ as a quotient space,  $\tau \colon \mW \to \mK$.    The essence  of the novel strategy behind the aperiodic property of $\Phi_t$  is perhaps best described by a quote from the paper by Matsumoto \cite{Matsumoto1995}:
\begin{quotation}
  We therefore must demolish the two closed orbits in the Wilson Plug beforehand. But producing a new plug will take us back to the starting line. The idea of Kuperberg is to \emph{let closed orbits demolish themselves}. We set up a trap within enemy lines and watch them settle their dispute while we take no active part. 
\end{quotation}

There are many choices made in the implementation of this strategy, where all such choices result in aperiodic flows. On the other hand, some of the choices appear to impact the  further dynamical properties of the Kuperberg flows, as will be discussed later.  We indicate in this section these alternate choices, and later formulate some of the questions which appear to be   important for further study.
Finally,  at the end of this section, we consider flows on plugs which violate the above strategy, where the traps for the periodic orbits are purposely not aligned. This results in what we call ``Derived from Kuperberg'' flows, or simply DK--flows for short, which were introduced in the work   \cite{HR2016b}.

\begin{figure}[!htbp]
\centering
{\includegraphics[width=72mm]{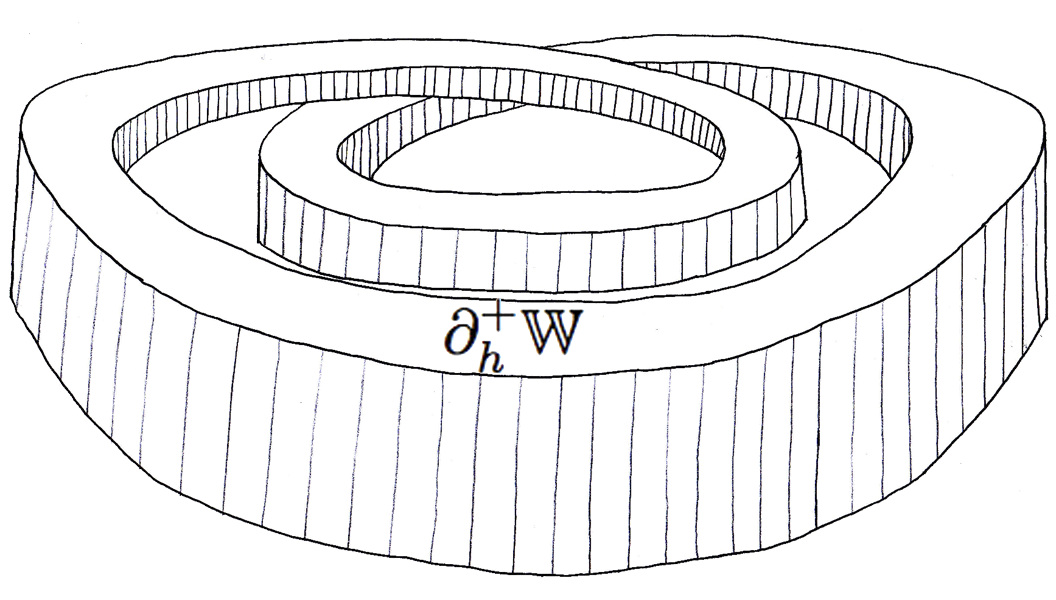}}
\caption{Embedding of Wilson Plug $\mW$ as a {\it folded figure-eight} \label{fig:8doblado} }
\vspace{-6pt}
\end{figure}

The construction of a Kuperberg Plug  $\mK$ begins with    the modified
Wilson Plug   $\mW$ with vector field $\cW$ constructed in Section~\ref{sec-wilson}.
The first step is to re-embed the manifold $\mW$   in $\mR^3$ as a {\it folded
 figure-eight}, as shown in Figure~\ref{fig:8doblado}, preserving the
vertical direction.

 The next step   is to   construct  two (partial)
insertions  of $\mW$ in itself,  so that each
periodic orbit of the Wilson flow is ``broken open'' by a trapped
orbit on the self-insertion.

\begin{figure}[!htbp]
\centering
{\includegraphics[width=50mm]{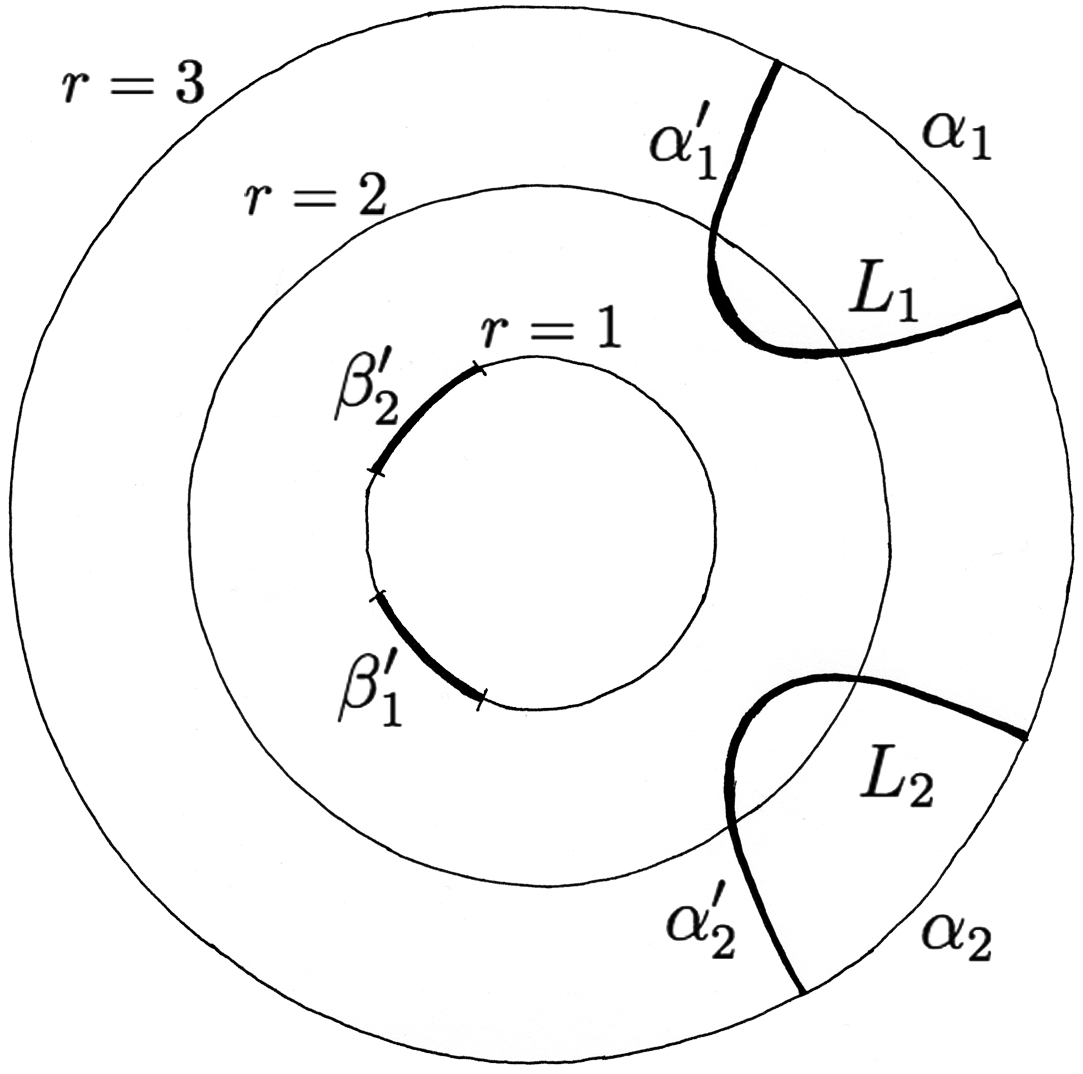}}
\caption{ The disks $L_1$ and $L_2$ \label{fig:insertiondisks}}
\vspace{-6pt}
\end{figure}

The construction begins with the choice    in the annulus $[1,3] \times \mS^1$  of  two   closed regions $L_i$, for $i=1,2$, which are topological disks. Each region has boundary defined by   two arcs: for $i=1,2$, $\alpha^\prime_i$ is the boundary contained in  the
interior of $[1,3] \times \mS^1$  and $\alpha_i$ in the outer boundary contained in the circle
$\{r=3\}$, as depicted in Figure~\ref{fig:insertiondisks}. In the work \cite{HR2016a} the choices for these arcs are given more precisely, though for our discussion here this is not necessary.

Consider the closed sets $D_i \equiv L_i \times[-2,2] \subset \mW$, for $i = 1,2$.  
Note that  each $D_i$ is homeomorphic to a closed $3$-ball, that $D_1 \cap D_2 = \emptyset$,  and each $D_i$ intersects the  cylinder $\cC = \{r=2\}$ in   a rectangle. Label   the top and bottom faces of these regions
\begin{equation}\label{eq-regions}
L_1^{\pm}=  L_1\times \{\pm 2\}  ~ , ~ 
L_2^{\pm} =    L_2\times \{\pm 2\}.
\end{equation}

The next step is to define  insertion maps   $\sigma_i  \colon D_i \to \mW $, for $i=1,2$, in such
a way that the periodic orbits $\cO_{1}$ and $\cO_{2}$ for the $\cW$-flow  intersect $\sigma_i(L_i^-)$ in   points
corresponding to   $\cW$-trapped points. 
Consider   two disjoint arcs $\beta_i'$ in the inner boundary circle
$\{r=1\}$ of $[1,3]\times \mS^1$, again as depicted in Figure~\ref{fig:insertiondisks}.

\begin{figure}[!htbp]
\centering
{\includegraphics[width=60mm]{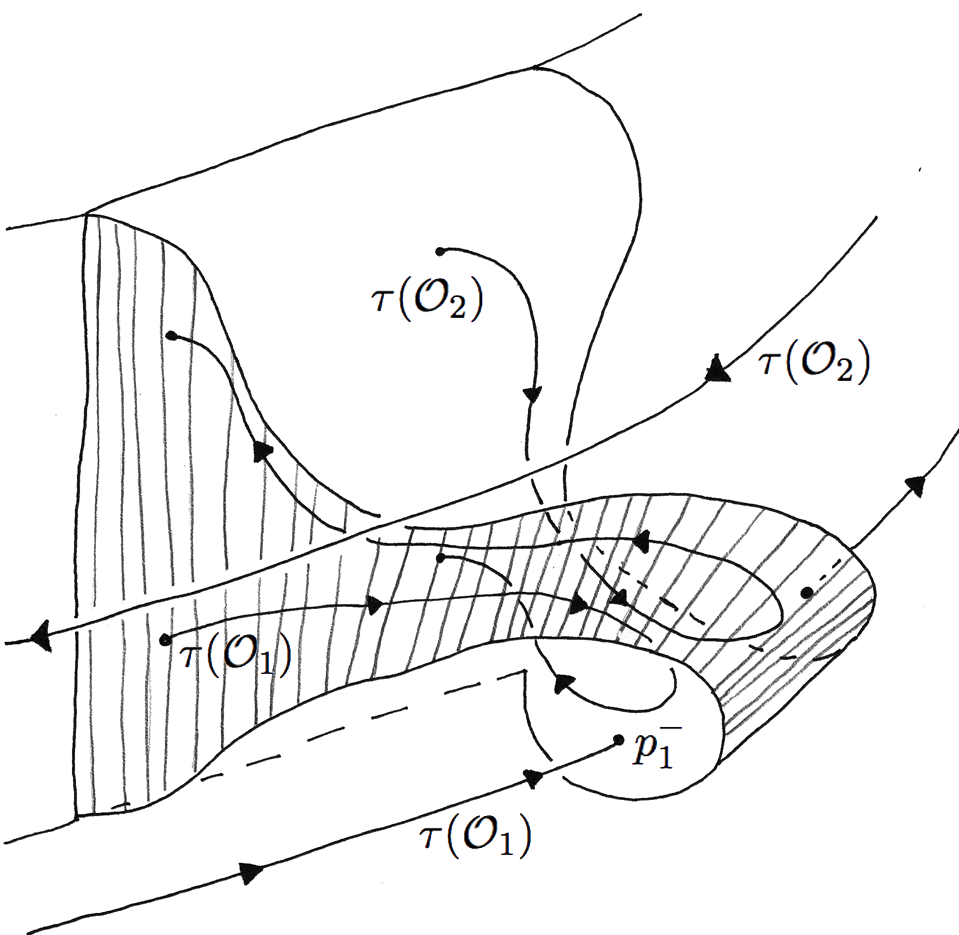}\hspace{20pt} \includegraphics[width=50mm]{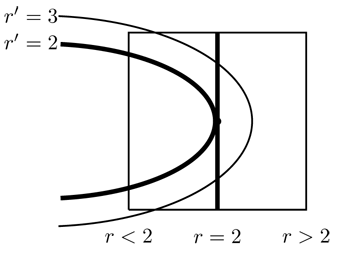}}
\caption{The image of $L_1\times [-2,2]$ under $\sigma_1$ and the radius function\label{fig:twisted} }
\vspace{-6pt}
\end{figure}

Now   choose a smooth family of
orientation preserving  diffeomorphisms $\sigma_i \colon \alpha_i'
\to \beta_i'$, $i=1,2$. 
Extend these maps to   smooth embeddings
$\sigma_i  \colon D_i \to  \mW $, for $i=1,2$, as illustrated on the left-hand-side of   
Figure~\ref{fig:twisted}. We require the following  conditions   for $i=1,2$:   
\begin{itemize}
\item[(K1)]   $\sigma_i(\alpha_i'\times z)=\beta_i'\times z$ for all $z\in [-2,2]$,   the interior arc $\alpha_i^\prime$ is mapped to a boundary arc $\beta_i'$. 
\item[(K2)]     $\cD_i = \sigma_i(D_i)$ then $\cD_1 \cap \cD_2 =\emptyset$;
\item[(K3)]    For every $x \in L_i$, the image   $\cI_{i,x} \equiv \sigma_i  (x \times
  [-2,2])$ is an arc contained in a trajectory of $\cW$;
\item[(K4)] $\sigma_1  (L_1 \times \{-2\}) \subset \{z < 0\}$  and $\sigma_2  (L_2 \times \{2\}) \subset \{z > 0\}$; 
\item[(K5)]  Each slice $\sigma_i  (L_i\times\{z\})$ is transverse to the vector field $\cW$, for all $-2\leq z \leq 2$. 
\item[(K6)]   $\cD_i$ intersects the periodic orbit $\cO_i$ and not $\cO_j$, for   $i \ne j$.
\end{itemize}
The ``horizontal faces'' of the embedded regions   $\cD_i \subset \mW$ are labeled by
\begin{equation}\label{eq-tongues}
\cL_1^{\pm}= \sigma_1(L_1^\pm)  ~ , ~   
\cL_2^{\pm} =   \sigma_2(L_2^\pm).
\end{equation}
Then the above assumptions imply that lower faces $\cL_1^{\pm}$ 
intersect  the first periodic orbit $\cO_1$ and are disjoint from the
second periodic orbit $\cO_2$, while the upper faces $\cL_2^{\pm}$
intersect   $\cO_2$ and are disjoint from  $\cO_1$.

The embeddings $\sigma_i$ are also required to satisfy two further conditions, which are the key to showing that the resulting Kuperberg flow   $\Phi_t$ is \emph{aperiodic}:   
\begin{itemize}
\item[(K7)] For $i=1,2$, the disk $L_i$ contains a point $(2,\theta_i)$ such that
  the image under $\sigma_i$ of the vertical segment
  $(2,\theta_i)\times[-2,2] \subset D_i \subset \mW$ is an arc
  $\{r=2\} \cap \{\theta_i^- \leq \theta \leq \theta_i^+\} \cap
  \{z=(-1)^i\}$  of the periodic orbit $\cO_i$.   
\item[(K8)] {\it Radius Inequality}: For all $x' = (r', \theta',  z') \in  L_i \times [-2,2]$, let $x = (r, \theta,z) = \sigma_i(r', \theta',  z') \in \cD_i$,  then $r < r'$ unless  $x' =  (2,\theta_i, z')$ and then $r=r'=2$.
\end{itemize}

The Radius Inequality (K8) is one of the most fundamental  concepts of Kuperberg's construction, and is illustrated   by the graph on the right-hand-side of in Figure~\ref{fig:twisted}.

Condition (K4) and the fact that  the flow of the vector field $\cW$ on $\mW$ preserves the radius coordinate on $\mW$, allow restating (K8) in the more concise form for points in the faces $\cL_i^-$ of the insertion regions $\cD_i$. For $x = (r, \theta, z) = \sigma_i(r', \theta',  z') \in \cD_i$ we have
 \begin{equation}\label{eq-radius}
r(\sigma_i^{-1}(x)) \geq r ~ {\rm for} ~ x \in \cL_i^- ~, ~ {\rm with }~ r(\sigma_i^{-1}(x)) = r ~{\rm if ~ and ~ only ~ if} ~ x = \sigma_i(2, \theta_i,  -2) ~. 
\end{equation}
The illustration of the radius inequality  in Figure~\ref{fig:twisted} is  an ``idealized'' case, as it implicitly assumes that the relation between the values of $r$ and $r'$ is ``quadratic'' in a neighborhood of the special points $(2,\theta_i)$, which is not required in order that  (K8) be satisfied.

\begin{figure}[!htbp]
\centering
{\includegraphics[width=120mm]{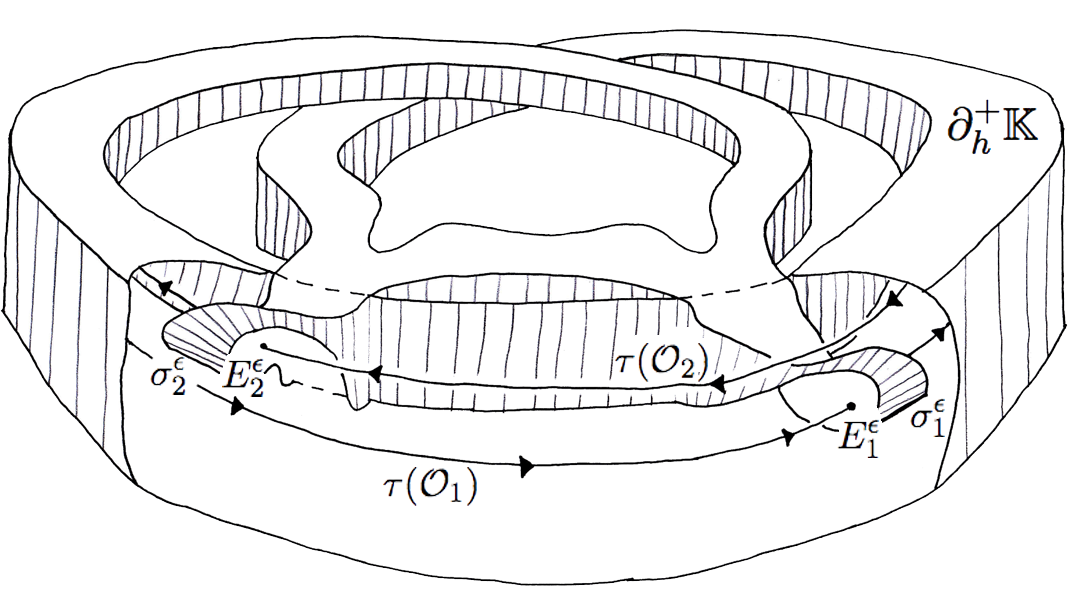}
 \caption{\label{fig:K} The Kuperberg Plug $\mK_\e$}}
\vspace{-6pt}
\end{figure}

Finally,  define $\mK$ to be the quotient manifold obtained from $\mW$ by
identifying the sets $D_i$ with $\cD_i$. That is, for each point $x
\in D_i$ identify $x$ with $\sigma_i(x) \in \mW$, for $i = 1,2$. 
 The restricted $\cW$-flow on the inserted disk $\cD_i = \sigma_i(D_i)$ is not compatible with the image of the restricted $\cW$-flow on $D_i$.  Thus, to obtain a smooth vector field $\cX$ from this construction,  it is necessary to modify $\cW$ on each insertion $\cD_i$. 
 The idea is to  replace  the vector field $\cW$  in the interior of each  region $\cD_i$ with the image vector field and smooth the resulting piecewise continuous flow \cite{Kuperberg1994,Ghys1995}.
Then the vector field $\cW^\prime$ on $\mW'$ descends to a smooth vector field on $\mK$ denoted by $\cK$, whose flow   is denoted by   $\Phi_t$. 
The family of \emph{Kuperberg Plugs} is the resulting space  $\mK \subset \mR^3$,   as   illustrated  in Figure~\ref{fig:K}. 

The images in $\mK$ of the cut-open periodic orbits from the Wilson flow $\Psi_t$ on $\mW$,  generate two orbits for the Kuperberg flow $\Phi_t$ on $\mK$, which are called  the \emph{special     orbits} for $\Phi_t$. These two special orbits play an absolutely central role in the study of the dynamics of the  flow $\Phi_t$. We   now state Kuperberg's main result:

\begin{thm}\cite{Kuperberg1994}\label{thm-aperiodic}
The flow $\Phi_t$ on $\mK$ satisfies the conditions on a plug, and has no periodic orbits.
\end{thm}

The papers   \cite{Ghys1995,Kuperberg1994} remark that a  Kuperberg
Plug  can also be constructed for which the manifold $\mK$ and its
flow $\cK$ are real analytic. An explicit construction of such a flow
is given in \cite[Section~6]{Kuperbergs1996}. There is the added
difficulty that the insertion of the plug in an analytic manifold must
also be analytic, which requires some subtlety. This is discussed in
\cite[Section~6]{Kuperbergs1996}, and also  in the second author's  Ph.D. Thesis \cite[Section~1.1.1]{Rechtman2009}.

Finally, we introduce a modification  to the above construction, for
which the periodic orbits of the Wilson flow are not necessarily broken open by
the trapped orbits of the inserted regions. 
 Let $\e$ be a fixed small constant, positive or negative. 
 Choose  smooth embeddings $\sigma_i^\e  \colon D_i \to  \mW $, for $i=1,2$, again as illustrated on the left-hand-side of   
Figure~\ref{fig:twisted}, which satisfy the conditions (K1) to (K6). In place of the conditions (K7) and (K8), we impose the modified conditions:   
\begin{itemize}
\item[(K7$\e$)] For $i=1,2$, the disk $L_i$ contains a point $(2,\theta_i)$ such that
  the image under $\sigma_i ^\e$ of the vertical segment
  $(2,\theta_i)\times[-2,2] \subset D_i \subset \mW$ is an arc of a
  $\cW$-orbit in 
  $\{r=2+\e\} \cap \{\theta_i^- \leq \theta \leq \theta_i^+\}$.   \\
\item[(K8$\e$)] {\it Parametrized Radius Inequality}: For all $x' = (r', \theta',  -2) \in  L_i^-$, let $x = (r, \theta,z) =
  \sigma_i^\e(r', \theta',  -2) \in \cL_i^{\e-}$,  then $r < r'+\e$ unless  $x' = (2,\theta_i, -2)$ and then $r=2+\e$.
\end{itemize}
 
\begin{figure}[!htbp]
\centering
\begin{subfigure}[c]{0.3\textwidth}{\includegraphics[height=32mm]{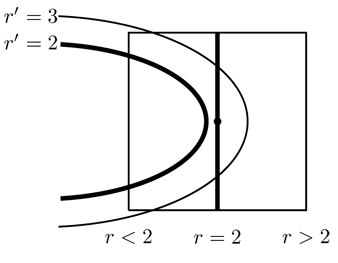}}\caption{$\e < 0$}\end{subfigure}
\begin{subfigure}[c]{0.3\textwidth}{\includegraphics[height=32mm]{pix/Figure07b.png}}\caption{$\e = 0$}\end{subfigure}
\begin{subfigure}[c]{0.3\textwidth}{\includegraphics[height=32mm]{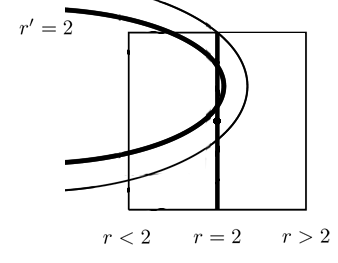}}\caption{$\e > 0$}\end{subfigure}
\caption{\label{fig:modifiedradius}  The modified radius inequality for  the cases $\e<0$, $\e=0$ and  $\e>0$}  
\vspace{-6pt}
\end{figure}

Observe that for $\e=0$, we recover the Radius Inequality (K8). Figure~\ref{fig:modifiedradius} represents the radius
inequality for the three cases where $\e<0$, $\e =0$,  and $\e>0$.

Again,   define $\mK_\e$ to be the quotient manifold obtained from $\mW$ by
identifying the sets $D_i$ with $\cD_i^\e = \sigma_i^\e(D_i)$.  
 Replace  the vector field $\cW$  on the interior of each  region $\cD_i^\e$ with the image vector field and smooth the resulting piecewise continuous flow, so that we obtain  a smooth vector field on $\mK_\e$ denoted by $\cK_\e$, whose flow   is denoted by   $\Phi_t^\e$. We say that  $\Phi_t^\e$ is a \emph{Derived from Kuperberg flow}, or a DK--flow.

 The dynamics of a DK--flow is actually quite simple in the case when $\e < 0$, as shown by the following result.

\begin{thm}\cite{HR2016b}\label{thm-periodic} 
Let $\e < 0$ and   $\Phi_t^\e$ be a DK--flow on $\mK_e$. Then it satisfies    the conditions on a plug, and moreover 
  the flow in the plug $\mK_\e$ has two periodic orbits that   bound an embedded invariant cylinder, and every other orbit belongs to the wandering set. 
\end{thm}
The proof of Theorem~\ref{thm-periodic} in \cite{HR2016b} uses the same  technical tools as developed in the previous works \cite{Kuperberg1994,Kuperbergs1996,Ghys1995,Matsumoto1995,HR2016a} for the study of the dynamics of Kuperberg flows.
In contrast,  the dynamics of a DK--flow  when $\e > 0$ can be quite
chaotic, having positive topological entropy and have an abundance of periodic orbits, as shown by the
construction of examples in \cite{HR2016b}.

  \section{Generic hypotheses}\label{sec-generic}

The construction of Kuperberg flows $\Phi_t$ on the plug $\mK$ which
satisfy Theorem~\ref{thm-aperiodic} involved multiple choices, which
do not change whether the resulting flows  are aperiodic, but do impact other dynamical properties of these flows.   In this section, we discuss in more detail these choices, and introduce the generic assumptions that were imposed in the works \cite{HR2016a,HR2016b}. The implications of these choices will be discussed in subsequent sections.
We first discuss the choices made in constructing the modified Wilson
plug, then consider the even wider range of choices involved with the
construction of the insertion maps.  Note that we discuss first the case for the traditional Kuperberg flows, and afterwards
discuss the variations for the case of the Derived from Kuperberg
flows (DK--flows).

 Recall that the   modified  Wilson vector field on $\mW$  is given in \eqref{eq-wilsonvector} by 
\begin{equation*} 
\cW =g(r, \theta, z)  \frac{\partial}{\partial  z} + f(r, \theta, z)  \frac{\partial}{\partial  \theta}
\end{equation*}
where the function $g(r, \theta, z) = g(r,  z)$ is  the suspension of the function $g \colon \bR \to [0,g_0]$   which is positive, except at the points $(2,  \pm 1)$, and symmetric about the line $\{z=0\}$. The function $f(r, \theta, z) = f(r, z)$ is assumed to satisfy the conditions (W1) to (W6), though conditions (W5) and (W6) are imposed   to simplify   calculations, and do not impact the aperiodic conclusion for the Kuperberg flows. 

Recall that  Figure~\ref{fig:flujocilin} illustrates the dynamics of
the flow   of $\cW$ restricted to the cylinders $\{r=cst.\}$ in $\mW$,
for various values of the radius. It is clear from these pictures that
the ``interesting'' part of the dynamics of this flow occurs on the
cylinders with radius near to $2$, and near the periodic orbits $\cO_i$ for $i=1,2$.

The points $(2, \pm 1) \in \bR$ are the local minima for the function $g$, and thus  its matrix of first derivatives must also vanish at these points, and the Hessian matrix of second derivatives must be positive semi-definite.
The generic property for such a function is that the Hessian matrix for $g$ at these points is positive definite. 
In the works \cite{HR2016a,HR2016b}, a more precise version of this was formulated:
 \begin{hyp}\label{hyp-genericW}
The function $g$ satisfies the following conditions:
\begin{equation}\label{eq-generic1}
g(r,z) = g_0 \quad \text{for} \quad  (r-2)^2 + (|z| -1)^2 \geq \e_0^2 
\end{equation}
where $0 < \e_0 < 1/4$ is sufficiently small. 
Moreover, we require that the Hessian   matrices of second partial derivatives for $g$ at the vanishing points $(2, \pm 1)$ are positive definite.
In addition, we   require   that $g(r,z)$ is monotone increasing as a function of the distance $\ds \sqrt{(r-2)^2 + (|z| -1)^2}$ from the   points $(2, \pm 1)$. 
\end{hyp}
 The conclusions of  Proposition~\ref{prop-wilsonproperties} do not require Hypothesis~\ref{hyp-genericW}, and so Theorem~\ref{thm-aperiodic} does not require it. On the other hand,  many of the results in  \cite{HR2016a,HR2016b} do require this generic hypothesis for their proofs, as it allows making estimates on the ``speed of ascent'' for the orbits of the Wilson flow near the periodic orbits.  

Hypothesis~\ref{hyp-genericW} implies a local quadratic estimate on the function $g$ near the points $(2, \pm 1)$ which is given as estimate  (94) in \cite{HR2016a}. We formulate   a more general version of this local estimate for $g$.

 \begin{hyp}\label{hyp-genericWn}
 Let $n \geq 2$ be an even integer. Assume there exists constants $\lambda_2 \geq \lambda_1> 0$   and $\e_0 > 0$ such that  
  \begin{equation}\label{eq-genericWn}
\lambda_1   \left( (r-2)^n + (|z|-1)^n \right)  \leq  g(r, z) \leq \lambda_2   \left( (r-2)^n + (|z|-1)^n \right)    \quad {\rm for} ~   \left( (r-2)^n + (|z|-1)^n \right)  \leq \e_0^{n+2} \ . 
\end{equation}
We then say that the resulting vector field $\cW$ on $\mW$ \emph{vanishes with order $n$}.
\end{hyp}
Hypothesis~\ref{hyp-genericW} implies that Hypothesis~\ref{hyp-genericWn} holds for $n=2$. This yields an estimate on the speed which the orbits of $\cW$ in $\mW$ for points with $z   \ne \pm 1$ approach the periodic orbits $\cO_i$ in forward or backward time, as discussed in detail in \cite[Chapter~17]{HR2016a}. When $n > 2$,  this speed of approach becomes slower and slower as $n$ gets larger. We can also allow for the case where $g$ has all partial derivatives vanishing at the points $(2,  \pm 1)$, in which case we say that the function $g$ vanishes to infinite order at the critical points, and we say that the resulting vector field $\cW$ on $\mW$ \emph{infinitely flat} at $\cO_i$ for $i=1,2$. In that case, the speed of approach of orbits of $\cW$ in $\mW$ become arbitrarily slow towards the periodic orbits.

The choices for the embeddings  $\sigma_i  \colon D_i \to  \mW $, for $i=1,2$, as illustrated on the left-hand-side of   
Figure~\ref{fig:twisted}, are more wide-ranging, and have  a fundamental influence  on the dynamics of the resulting Kuperberg flows on the quotient space $\mK$. 
 We first impose a ``normal form'' condition on the insertions, which does not have significant impact on the dynamics, but allows a more straightforward formulation of the other properties of the insertion maps. 

 Let $(r, \theta,z)  = \sigma_i(x') \in \cD_i$ for    $i=1,2$, where     $x' = (r', \theta',  z') \in  D_i$ is a point in the domain of $\sigma_i$.
Let $\pi_z(r, \theta, z) = (r, \theta,-2)$ denote the projection of $\mW$ along the $z$-coordinate. 
We assume that $\sigma_i$ restricted to the bottom face, $\sigma_i \colon L_i^- \to \mW$,  has image  transverse to the vertical fibers of $\pi_z$. This normal form can be achieved by an isotopy of a given embedding along the flow lines of the vector field $\cW$, so does not change the orbit structure of the resulting vector field on the plug $\mK$.

 The above transversality assumption implies that   $\pi_z \circ \sigma_i \colon L_i^- \to \mW$ is a diffeomorphism into the face $\partial_h^- \mW$, with   image  denoted by $\fD_i \subset \partial_h^- \mW$. 
Then let $\ds \vartheta_i = (\pi_z \circ \sigma_i)^{-1} \colon \fD_i \to L_i^-$ denote the inverse map, so we have: 
\begin{equation}\label{eq-coordinatesVT}
\vartheta_i(r, \theta,-2) = (r(\vartheta_i(r,\theta,-2)), \theta(\vartheta_i (r,\theta,-2)), -2)  = (R_{i,r}(\theta)  , \Theta_{i,r}(\theta), -2)  ~.
\end{equation}
 We can then formalize in terms of the maps $\vartheta_i$ the   assumptions on   the insertion maps $\sigma_i$ that are intuitively 
implicit in  Figure~\ref{fig:twisted}, and will be assumed for all insertion maps considered.

 \begin{hyp} [{\it Strong Radius Inequality}]  \label{hyp-SRI}
  For   $i=1,2$, assume that: 
  \begin{enumerate}
\item \label{item-SRI-2} $\sigma_i \colon L_i^- \to \mW$ is transverse to the fibers of $\pi_z$; 
\item \label{item-SRI-1} $r = r(\sigma_i(r',\theta',z'))< r'$, except for  $(2,\theta_i, z')$ and  then  $z(\sigma_i(2,\theta_i, z')) = (-1)^i$; 
  \item \label{item-SRI-3} $\Theta_{i,r}(\theta) = \theta(\vartheta_i(r, \theta,-2))$ is an increasing function of $\theta$ for each fixed $r$; 
\item  \label{item-SRI-4} $R_{i,r}(\theta) = r(\vartheta_i(r,  \theta,-2))$ has non-vanishing derivative for $r=2$, except for the case of  $\overline{\theta_i}$ 
  defined by  $\vartheta_i(2,\overline{\theta_i},-2)= (2,\theta_i,-2)$; 
  \item For $r$ sufficiently close to $2$, we require that the $\theta$ derivative of $R_{i,r}(\theta)$ vanish at a unique point denoted by $\overline{\theta}(i,r)$.
\end{enumerate}
Consequently, each  surface $\cL_i^-$   is transverse to the coordinate vector fields $\partial/\partial \theta$ and $\partial/\partial z$ on $\mW$.
  \end{hyp}

 The illustration of the image of the curves $r'=2$ and $r'=3$  on the right-hand-side of  Figure~\ref{fig:twisted} suggests that these curves have ``parabolic shape''.  We formulate this notion more precisely using the function $\vartheta_i(r, \theta,-2)$ defined by \eqref{eq-coordinatesVT}, and introduce the more general hypotheses they may satisfy.
  Recall that $\e_0 > 0$ was introduced in Hypothesis~\ref{hyp-genericW}.

\begin{hyp}\label{hyp-polynomialn}   
Let $n \geq 2$ be an even integer.
For  $i = 1,2$,   $2 \leq r_0 \leq 2 + \e_0$ and $\theta_i -\e_0 \leq \theta \leq \theta_i + \e_0$, assume that   
\begin{equation}\label{eq-polynomialn}
\frac{d}{d\theta} \Theta_{i,r_0}(\theta) > 0 \quad, \quad  \frac{d^n}{d\theta^n} R_{i,r_0}(\theta) > 0  \quad, \quad \frac{d}{d\theta} R_{i,r_0}(\overline{\theta_i}) = 0 \quad, \quad \frac{d^{\ell}}{d\theta^{\ell}}  R_{i,r_0}(\overline{\theta_i}) = 0 ~ {\rm for} ~ 1 \leq \ell < n \ .
\end{equation}
where $\overline{\theta_i}$ satisfies $\vartheta_i(2,\overline{\theta_i},-2)= (2,\theta_i,-2)$. Thus for $2 \leq r_0 \leq 2+ \e_0$,  the graph of $R_{i,r_0}(\theta)$ is   convex upwards with vertex    at   $\theta = \overline{\theta_i}$. 
\end{hyp}
    
 In the case where $n=2$,  Hypothesis~\ref{hyp-polynomialn}  implies
 that all of the level curves $r'=c$, for $2 \leq c \leq 2 + \e_0$,
 have parabolic shape, as the illustration in
 Figure~\ref{fig:twisted} suggests.   On the other hand, for $n > 2$
 the level curves $r'=c$ have higher order contact with the vertical
 lines of constant radius in Figure~\ref{fig:twisted}, and in this case, many of the dynamical properties of the resulting flow $\Phi_t$ on $\mK$ are not well-understood. 
 
We can now define what is called a generic Kuperberg flow in the work \cite{HR2016a}.
 \begin{defn}\label{def-generic}
 A Kuperberg flow $\Phi_t$ is \emph{generic}   if the Wilson flow $\cW$ used in the construction of the vector field $\cK$   satisfies Hypothesis~\ref{hyp-genericW}, and the insertion maps $\sigma_i$ for $i=1,2$ used in the construction of $\mK$ satisfies 
 Hypotheses~\ref{hyp-SRI},  and Hypotheses~\ref{hyp-polynomialn} for $n=2$. That is,  the singularities for the vanishing of the vertical component $g \cdot   \partial/ \partial  z$ of the vector field $\cW$ are of quadratic type, and the insertion maps used to construct $\mK$ yield  quadratic  radius functions near the special points.
 \end{defn}
 
 \eject

Recall that the insertion maps  for a Derived from Kuperberg flow as introduced in Section~\ref{sec-kuperberg} are denoted by 
$\sigma_i^\e  \colon D_i \to  \mW $, for $i=1,2$.   It is assumed that these maps satisfy   the modified conditions    (K7$\e$) and (K8$\e$).  The illustrations of the radius inequality in Figure~\ref{fig:modifiedradius} again suggest that the images of the curves $r'=c$ are of ``quadratic type'', though the vertex of the image curves need no longer    be at a special point. 
We again assume the    insertion maps 
 $\sigma_i^\e  \colon L_i^- \to \mW$ are transverse to the fibers of the projection map  $\pi_z \colon \mW \to \partial_h^-\mW$  along the $z'$-coordinate. 
Then we can define  the inverse map $\vartheta_i^\e = (\pi_z \circ \sigma_i^\e)^{-1} \colon \fD_i \to L_i^-$ 
and   express the inverse map $x' = \vartheta_i^\e(x)$ in polar coordinates as:
\begin{equation}\label{eq-coordinatesVTe}
x' = (r',\theta', -2) = \vartheta_i^\e(r, \theta,-2) = (r(\vartheta_i^\e(r,\theta,-2)), \theta(\vartheta_i^\e (r,\theta,-2)), -2)  = (R_{i,r}^\e(\theta)  , \Theta_{i,r}^\e(\theta), -2)  ~.
\end{equation}
  Then the level curves $r' = c$ pictured in Figure~\ref{fig:modifiedradius} are given by the maps $\theta' \mapsto \pi_z(\sigma_i^\e(c, \theta',  -2)) \in \partial_h^- \mW$.

We note    a straightforward consequence of the Parametrized Radius Inequality (K8$\e$). 
Recall that $\theta_i$ is the radian coordinate specified in (K8$\e$) such that for $x' = (2, \theta_i,  -2) \in  L_i^-$ we have 
$r(\sigma_i^\e(2, \theta_i,  -2)) = 2+\e$.

\begin{lemma}\cite[Lemma~6.1]{HR2016b}\label{lem-r0}
For $\e>0$ there exists $2+\e <r_\e<3$ such that $r(\sigma_i^\e(r_\e,\theta_i, -2))= r_\e$.
\end{lemma}

 We then add an additional assumption on the   insertion maps $\sigma_i^\e$ for $i=1,2$ which specifies the qualitative behavior of the radius function for $r \geq r_\e$.  
 \begin{hyp} \label{hyp-monotone}
If $r_\e$ is the smallest  $2+\e <r_\e<3$  such that
$r(\sigma_i^\e(r_\e,\theta_i, -2))= r_\e$. Assume that $r(\sigma_i^\e(r,\theta_i, -2)) < r$ for $r > r_\e$.  
 \end{hyp}
 The conclusion of Hypothesis~\ref{hyp-monotone}    is implied by the Radius Inequality for the case $\e=0$, but does not follow from the condition (K8$\e$) when $\e > 0$. It is imposed to eliminate some of the possible pathologies in the behavior of the orbits of the DK--flows.

 We can now formulate the analog for DK--flows of the Hypothesis~\ref{hyp-SRI}, which imposes      uniform conditions on the derivatives of the maps $\vartheta_i^\e$.   
 Recall that $0 < \e_0 < 1/4$ was specified in Hypothesis~\ref{hyp-genericW}, and we assume that $0 < \e < \e_0$.

 \begin{hyp} [{\it Strong Radius Inequality}]  \label{hyp-SRIe}
  For   $i=1,2$, assume that: 
  \begin{enumerate}
\item \label{item-SRI-2e} $\sigma_i^\e \colon L_i^- \to \mW$ is transverse to the fibers of $\pi_z$; 
\item \label{item-SRI-1e} $r = r(\sigma_i^\e(r',\theta',z))< r+\e$, except for  $x' = (2,\theta_i, z)$ and  then  $r=2+\e$; 
  \item \label{item-SRI-3e} $\Theta_{i,r}^\e(\theta)$ is an increasing function of $\theta$ for each fixed $r$; 
\item  \label{item-SRI-4e}  For $2-\e_0 \leq r \leq 2+\e_0$ and  $i = 1,2$, assume that   $R_{i,r}^\e(\theta)$ has non-vanishing derivative, except when    $\theta = \overline{\theta_i}$ as 
  defined by  $\vartheta_i^\e(2+\e,\otheta_i,-2)= (2,\theta_i,-2)$;   
  \item For $r$ sufficiently close to  $2+\e$, we require that the $\theta$ derivative of $R_{i,r}^\e(\theta)$ vanishes at a unique point denoted by $\overline{\theta}(i,r)$.
\end{enumerate}
  \end{hyp}

Note that  Hypotheses~\ref{hyp-monotone} and \ref{hyp-SRIe} combined  imply  that $r_\e$ is the unique value of $2+\e <r_\e<3$ for which  $r(\sigma_i^\e(r_\e,\theta_i, -2))= r_\e$. 
 We can then formulate the   analog of Hypothesis~\ref{hyp-polynomialn}.
 
\begin{hyp}\label{hyp-polynomialne}   
Let $n \geq 2$ be an even integer.
For    $2 - \e_0 \leq r_0 \leq 2 + \e_0$ and $\theta_i -\e_0 \leq \theta \leq \theta_i + \e_0$, assume that   
\begin{equation}\label{eq-polynomialne}
\frac{d}{d\theta} \Theta_{i,r_0}^\e(\theta) > 0 \quad, \quad  \frac{d^n}{d\theta^n} R_{i,r_0}^\e(\theta) > 0  \quad, \quad \frac{d}{d\theta} R_{i,r_0}^\e(\overline{\theta_i}) = 0 \quad, \quad \frac{d^{\ell}}{d\theta^{\ell}}  R_{i,r_0}^\e(\overline{\theta_i}) = 0 ~ {\rm for} ~ 1 \leq \ell < n \ .
\end{equation}
where $\overline{\theta_i}$ satisfies $\vartheta_i^\e(2,\overline{\theta_i},-2)= (2,\theta_i,-2)$. Thus for $2 - \e_0 \leq r_0 \leq 2+ \e_0$,  the graph of $R_{i,r_0}^\e(\theta)$ is   convex upwards with vertex    at   $\theta = \overline{\theta_i}$. 
\end{hyp}
  
Finally, we have the definition of the generic DK--flows studied in \cite{HR2016b}.
 \begin{defn}\label{def-genericDK}
 A DK--flow $\Phi_t^\e$ is \emph{generic}   if the Wilson flow $\cW$ used in the construction of the vector field $\cK$   satisfies Hypothesis~\ref{hyp-genericW}, and the insertion maps $\sigma_i^\e$ for $i=1,2$ used in the construction of $\mK$ satisfies 
 Hypotheses~\ref{hyp-SRIe},  and Hypotheses~\ref{hyp-polynomialne} for $n=2$. 
  \end{defn}

 \eject
 
  \section{Wandering and minimal sets}\label{sec-minimalset}
 
We next   discuss some of the basic topological dynamics properties the   Kuperberg flows in plugs.
Our main interest is in the   asymptotic behavior of their orbits, especially the non-wandering   and   wandering sets for the flow.   
There is an additional subtlety in these considerations, in that many orbits for the flow in a plug may escape from the plug, while other orbits are trapped in either the forward or backward directions, or possibly both.  We also  recall the results about the uniqueness of the minimal set.  First we recall some of the basic concepts for the flow in a plug.

Recall that   $\cD_i = \sigma_i(D_i)$ for $i =1,2$   are   solid $3$-disks embedded in $\mW$.
 Introduce  the sets:
\begin{equation}\label{eq-notchedW}
\mW' ~   \equiv   ~ \mW - \left \{ \cD_1  \cup \cD_2 \right\} \quad , \quad 
\wmW~   \equiv   ~ \overline{\mW - \left \{ \cD_1  \cup \cD_2 \right\}} ~ .
\end{equation}
The compact space $\wmW \subset \mW$ is the result of  ``drilling out''  the interiors of $\cD_1$ and $ \cD_2$.

  Let   $\tau \colon \mW  \to \mK$ denote the quotient map. Note that the  restriction $\tau' \colon \mW' \to \mK$ is injective and onto, while  for $i = 1,2$,  the map $\tau$   identifies   a point $x \in D_i$ with its image $\sigma_i(x)   \in \cD_i$.
    Let   $(\tau')^{-1}  \colon \mK \to \mW'$ denote the inverse map, which followed  by the inclusion $\mW' \subset \mW$, yields  the (discontinuous) map  $\tau^{-1} \colon \mK \to \mW$, where  $i=1,2$, we have:
\begin{equation}\label{eq-radiusdef}
\tau^{-1}(\tau(x)) =x ~ {\rm for} ~ x \in D_i  ~ ,~ {\rm and} ~   \sigma_i(\tau^{-1}(\tau(x))) = x ~ {\rm for} ~ x \in \cD_i~.
\end{equation}

Consider the embedded disks  $\cL_i^{\pm} \subset \mW$ defined by \eqref{eq-tongues}, which appear  as the   faces of the
insertions  in $\mW$.  Their images in the quotient manifold $\mK$ are denoted by:
\begin{equation}\label{eq-sections}
E_1= \tau(\cL_1^-) ~ , ~ S_1= \tau(\cL_1^+) ~ , ~   
E_2= \tau(\cL_2^-) ~ , ~ S_2= \tau(\cL_2^+)   ~ .
\end{equation}
Note that $\tau^{-1}(E_i) = L_i^-$, while $\tau^{-1}(S_i) = L_i^+$.  

  The {\it transition   points}  of an orbit of $\Phi_t$ are those points where the orbit intersects
 one of  the sets $E_i$ or $S_i$ for  $i=1,2$,  or is contained in a boundary component  $\partial_h^-\mK$ or $\partial_h^+\mK$.
 The transition points are classified as either \emph{primary} or \emph{secondary}, where $x\in \mK$ is: 
\begin{itemize}
\item a \emph{primary entry  point} if $x\in \partial_h^-\mK$;
\item  a \emph{primary exit  point} if $x \in \partial_h^+\mK$;
\item a \emph{secondary entry  point} if $x \in E_1 \cup E_2$;
\item a \emph{secondary exit  point}  $x \in S_1 \cup S_2$.
\end{itemize}
If a $\Phi_t$-orbit of a point $x \in \mK$ contains no transition points, then the restriction $\tau^{-1}(\Phi_t(x))$ is a continuous function of $t$, and in fact is contained in the $\Psi_t$-orbit of $x' = \tau^{-1}(x) \in \mW$.

Recall that $r \colon \mW \to [1,3]$ is the radius coordinate on $\mW$. Define the   (discontinuous)   radius  coordinate    $r \colon \mK \to [1,3]$,  where for $x \in \mK$ set $r(x) = r(\tau^{-1}(x))$. 
Then for $x \in \mK$  set     $\rho_x(t) \equiv r(\Phi_t(x))$, which is   the radius coordinate function along the $\cK$-orbit of $x$.   Note that if $\Phi_{t_0}(x)$ is not an entry/exit point, then the    function $\rho_x(t)$  is locally constant at $t_0$. On the other hand,  if $t_0$ is a point of discontinuity for  $\Phi_{t}(x)$, then $x_0 = \Phi_{t_0}(x)$  must be a secondary entry or exit point.

These properties of the radius function along orbits of the flow $\Phi_t$ gives a strategy for the study of the dynamics of the flow, and in fact provides the key technique in \cite{Kuperberg1994} used  to prove that the flow is aperiodic. A key idea is to index the points along the orbit of a point $x \in \mK$ by the intersections with the sets $E_1 \cup E_2$, for which the index increases by $+1$, or by their intersection with the sets $S_1 \cup S_2$, for which the index decreases by $-1$. This yields the integer-valued level function $n_x(t)$ which has $n_x(0) = 0$.

      Recall that  $\cO_i$ for $i=1,2$ denotes the periodic orbits for the Wilson flow on $\mW$, so that each intersection $\cO_i \cap \mW'$ consists of an open connected arc with endpoints $\cL_i^{\pm} \cap \cO_i$.
The  \emph{special entry/exit points} for the flow $\Phi_t$ are the images, for $i=1,2$, 
\begin{equation}\label{eq-special}
p_i^{-} = \tau(\cO_i \cap \cL_i^{-}) \in E_i ~ , ~ p_i^{+} = \tau(\cO_i \cap \cL_i^{+}) \in S_i ~ .
\end{equation}
Note that by definitions and the Radius Inequality, we have $r(p_i^{\pm}) = 2$ for $i=1,2$. 

 We now recall the results for the minimal set of Kuperberg flows based on the combined results from the works  \cite{Ghys1995,Kuperberg1994,Kuperbergs1996,Matsumoto1995}.
 It was observed by Kuperberg in \cite{Kuperberg1994}  that for $x \in
 \mK$ with $r(x) =2$, then either its forward orbit $\{\Phi_t(x) \mid
 t\geq 0 \}$ contains a special point in its closure, or this is true
 for the backward orbit $\{\Phi_t(x) \mid t\leq 0 \}$, or both
 conditions hold. Also, for $x \in \mK$ if the radius function
 $\rho_x(t) \geq c$ for some $c > 2$, then the orbit of $x$ escapes in
 finite time in both forward and backward directions.  It follows from this that for $x
 \in \mK$ with $r(x) > 2$ and whose orbit is infinite in either forward or
 backward directions, then its orbit closure  must contain  at least
 one of the special orbits. 
 
   It was observed in Matsumoto \cite{Matsumoto1995} that there is an
   open set of primary entry points with radius less than $2$ whose
   forward orbits are non-recurrent and yet accumulate on the special
   orbits. Ghys showed in \cite[Th\'eor\`eme, page 301]{Ghys1995} that
   if $x \in \mK$  does not escape from $\mK$ in a finite time, either
   forward or backward, then the orbit of the point accumulates on (at
   least one) of the special orbits.
These results combined imply that       a Kuperberg flow has   a unique minimal  set contained in the interior of $\mK$.

We state these results more succinctly as follows. 
  Define the following orbit closures in $\mK$:
\begin{equation}  \label{eq-minset1}
\Sigma_1 ~   \equiv   ~ \overline{ \{\Phi_t(p_1^-) \mid -\infty < t < \infty \}}     \quad , \quad 
\Sigma_2 ~   \equiv   ~ \overline{ \{\Phi_t(p_2^-) \mid -\infty < t < \infty \}}     ~ .
\end{equation}

\begin{thm} \cite[Theorem~8.2]{HR2016a} \label{thm-minimal} 
For the closed   sets   $\Sigma_i$ for $i=1,2$ we have:
\begin{enumerate}
\item $\Sigma_i$ is $\Phi_t$-invariant;
\item $r(x) \geq 2$ for all $x \in \Sigma_i$;
\item $\Sigma_1 = \Sigma_2$ and we set  $\Sigma = \Sigma_1 = \Sigma_2$;
\item Let  $\cZ \subset \mK$ be a closed invariant set for $\Phi_t$ contained in the interior of $\mK$, then $\Sigma \subset \cZ$;
\item $\Sigma$ is the unique minimal set for $\Phi_t$.
\end{enumerate}
\end{thm}

Note that for the Wilson plug   \cite{Wilson1966}, the flow has two
minimal sets, consisting of closed orbits, while the   Schweitzer
plug    \cite{Schweitzer1974} has also two minimal sets homeomorphic
to a Denjoy minimal set in the $2$-torus. On the other hand, the topological type of the unique minimal set $\Sigma$ for a Kuperberg flow is extraordinarily complicated, and seems to require additional generic hypotheses on the construction of the flow to gain a deeper understanding of its topological properties.

The orbits of the Kuperberg flow are divided into those which are finite, forward or backward trapped, or trapped in both directions and so infinite. A point $x \in \mK$  is \emph{forward wandering} if there exists an open set $x \in U \subset \mK$ and $T_U > 0$ so that for all $t \geq T_U$ we have $\Phi_t(U) \cap U = \emptyset$. Similarly, $x$ is \emph{backward wandering} if there exists an open set $x \in U \subset \mK$ and $T_U < 0$ so that for all $t \leq T_U$ we have $\Phi_t(U) \cap U = \emptyset$. A point $x$ with infinite orbit   is \emph{wandering} if it is forward and backward wandering. 
Define the following subsets of $\mK$:
\begin{eqnarray*}
\fW^{0} & \equiv & \{ x \in \mK \mid x ~ {\rm orbit ~ is ~ finite}\} \label{eq-wanderingfinite}\\
\fW^+ & \equiv & \{ x \in \mK \mid x ~ {\rm orbit ~ is ~ forward ~   wandering}\} \label{eq-wandering+}\\
\fW^- & \equiv & \{ x \in \mK \mid x ~ {\rm orbit ~ is ~ backward ~   wandering}\} \label{eq-wandering-}\\
\fW^{\infty} & \equiv & \{ x \in \mK \mid x ~ \mbox{ is forward and backward wandering} \} \label{eq-wanderinginfinite}
\end{eqnarray*}
Note that $x \in \fW^{0}$ if and only if the orbit of $x$ escapes through $\partial_h^+ \mK$ in forward time, and escapes though $\partial_h^- \mK$ in backward time. 
Define 
\begin{equation}\label{eq-dynamicdecomp}
\fW    ~ = ~  \fW^{0}    \cup \fW^+ \cup \fW^- \cup \fW^{\infty} \quad ; \quad  \Omega = \mK -  \fW .
\end{equation}

 The set $\Omega$ is called the \emph{non-wandering} set for $\Phi_t$, is    closed and $\Phi_t$-invariant. A point $x$ with forward trapped orbit    is characterized by the property: 
  $x \in \Omega$  if   for all $\e > 0$ and $T > 0$, there exists $y$
  and $t > T$ such that $d_{\mK}(x,y) < \e$ and $d_{\mK}(x, \Phi_t(y))
  < \e$, where $d_{\mK}$ is a distance function on $\mK$.
There are obvious corresponding statements for points which are   backward trapped or infinite. 
Here are some of the properties of the wandering and non-wandering sets for Kuperberg flows. 
The proofs can be found in  \cite[Chapter~8]{HR2016a}. 
  
\begin{lemma} \label{lem-wanderingboundary}
 If  $x \in \mK$   is   a primary entry or exit point, then $x \in
 \fW^+$ or $\fW^-$.  
\end{lemma}
 
\begin{lemma} \label{lem-wanderinginfinite}
For each $x \in \Omega$, the $\Phi_t$-orbit of $x$  is infinite.
\end{lemma}
 
 \begin{prop} \label{prop-wanderinginterior}
$\Sigma \subset \Omega \subset \{x \in \mK \mid r(x) \geq 2\}$. 
\end{prop}

Finally, let us recall a result of Matsumoto:
\begin{thm}\cite[Theorem~7.1(b)]{Matsumoto1995}
The sets $\fW^{\pm}$ contain interior points.
\end{thm}
This implies the following important consequence:
\begin{cor}
The flow $\Phi_t$ cannot preserve any smooth invariant measure on $\mK$ which assigns positive mass to any neighborhood of  a special point. 
\end{cor}

At this point it is worth mentioning one of the big open problems
regarding the existence of periodic orbits for flows on
3-manifolds. For volume preserving vectors fields that are at least
$C^2$, it is unknown whether they must have periodic
orbits. G. Kuperberg proved two important results: first every
3-manifold admits a volume-preserving $C^\infty$ vector field that
has a finite number of periodic orbits, and second that there are volume-preserving $C^1$-flows without periodic orbits \cite{KuperbergG1996}. These results
are based on the use of plugs.

  \section{Zippered laminations}\label{sec-lamination}

We next introduce the $\Phi_t$-invariant embedded
 surface $\fM_0$ and its closure $\fM$,  and  discuss the relation
 between the minimal set $\Sigma$   and the space $\fM$.   The
 existence of this compact  connected subset $\fM$ which is invariant for the
 Kuperberg flow $\Phi_t$ is a remarkable consequence  of the
 construction, and is   the key to a deeper understanding of the properties of the minimal set $\Sigma$ of $\Phi_t$. We then 
  give an overview of
 the structure theory for $\fM_0$ which plays a fundamental role in analyzing the dynamical properties of Kuperberg flows.

Recall that the Reeb cylinder   $\cR  \subset \mW$ is bounded by the
two periodic orbits $\cO_1$ and $\cO_2$ for the Wilson flow $\Psi_t$
on $\mW$. The cylinder $\cR$ is itself invariant under this  flow, and
for a point $x\in \mW$ with $r(x)$ close to $2$, the $\Psi_t$-orbit of
$x$   has increasing long orbit segments which shadow the periodic
orbits. Since the special orbits in $\mK$ contain the intersection
$\tau(\cO_i\cap\mW')$ for $i=1,2$, they are contained in the $\cK$-orbits of the
set obtained by flowing the Reeb cylinder $\tau(\cR\cap \mW')$.

\begin{figure}[!htbp]
\centering
{\includegraphics[width=60mm]{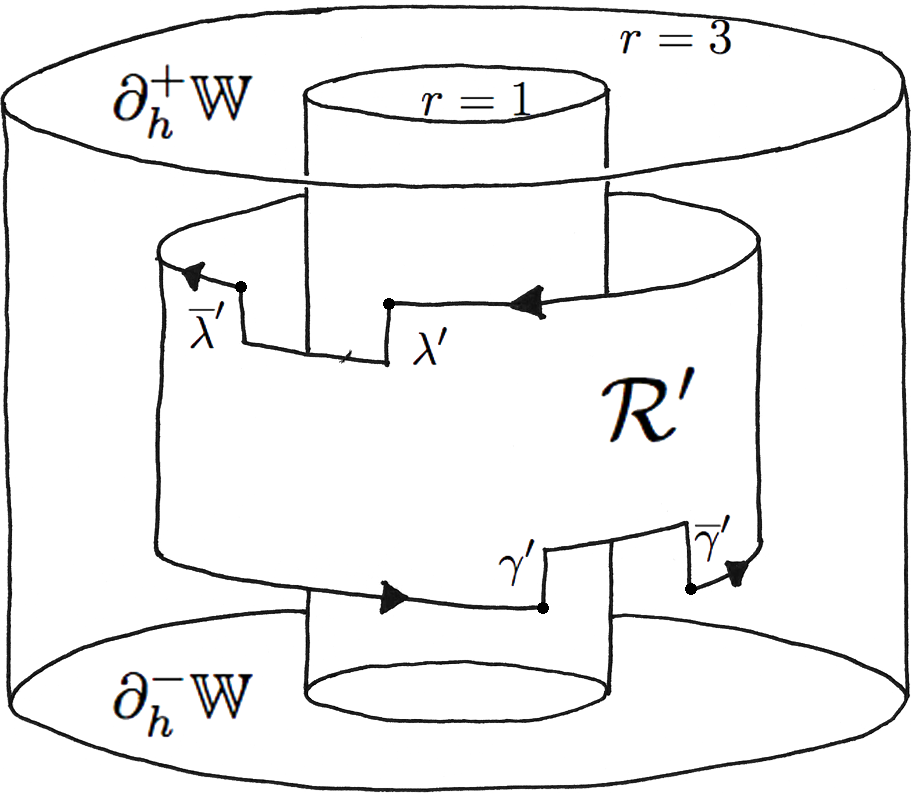}}
\caption{\label{fig:notches} The notched cylinder $\cR'$ embedded in $\mW$}
\vspace{-6pt}
\end{figure}

Introduce the 
 \emph{notched Reeb cylinder},  $\cR' = \cR \cap \mW'$, which has two
closed ``notches''   removed from   $\cR$ where it intersects the
closed insertions $\cD_i \subset \mW$ for $i=1,2$.
Figure~\ref{fig:notches}  illustrates the cylinder $\cR'$ inside $\mW$. 
The boundary segments  $\gamma'$ and $\lambda'$ labeled in Figure~\ref{fig:notches} satisfy $\gamma' \subset \cL_1^-$ and $\lambda' \subset \cL_2^-$, while   the boundary segments  $\ovg'$ and $\ovl'$ labeled in Figure~\ref{fig:notches} satisfy $\ovg' \subset \cL_1^+$ and $\ovl' \subset \cL_2^+$.  A basic observation is that these curves are each transverse to the restriction of the Wilson flow to the cylinder $\cR$.

The map $\tau \colon \cR' \to \mK$ is an embedding, so the 
  $\Phi_t$-flow of  $\tau(\cR') \subset \mK$  is   an   embedded surface,  
 \begin{equation}\label{eq-lamination}
\fM_0 ~ \equiv ~  \{\Phi_t(\tau(\cR')) \mid -\infty < t < \infty \}   ~.
\end{equation}
 The ``boundary'' of $\fM_0$ consists of the two   special orbits   in $\mK$ obtained by the $\Phi_t$-flows of the arcs $\tau(\cO_i \cap \mW')$ for $i=1,2$, so that $\fM_0$ is an ``infinite  bordism'' between the two special orbits of the flow $\Phi_t$. 
  Thus, the closure  $\fM = \overline{\fM_0}$ is a flow invariant,
  compact connected subset of $\mK$, which contains the closure of the
  special orbits, hence by Theorem~\ref{thm-minimal},   the minimal
  set $\Sigma \subset \fM$.  A fundamental problem is then to give a
  description of the topology and geometry of the space $\fM$.  
  The question of when $\Sigma=\fM$
  is treated in Section~\ref{sec-denjoy}, while in this section we concentrate on the properties
  of $\fM$.
  
  The key to understanding the structure of the space $\fM$ is to
  analyze  the structure  of $\fM_0$ and its embedding in  $\mK$. This
  analysis    is based on a simple observation, that the images
  $\tau(\gamma'), \tau(\lambda') \subset \fM$ are curves transverse to
  the flow $\Phi_t$ and contained in the region $\{x \in \fM \mid r(x)
  \geq 2\}$. Moreover, for a point $x \in \tau(\gamma')$   with $r(x)
  > 2$, there is a finite $t_x > 0$ such that    $\Phi_{t_x}(x) \in
  \tau(\ovg')$. That is, the flow across the notch in $\tau(\cR')$
  with boundary  curve $\tau(\gamma')$ closes up by returning to the
  facing boundary curve $\tau(\ovg')$, \emph{unless $r(x)=2$ and then
    $x$ is the special point $p_1^-$}. A similar remark holds for the notch in $\cR'$ with boundary curves $\lambda', \ovl'$.    
  It follows from the proof of the above remarks that  we can use a recursive approach to analyze the submanifold $\fM_0$,  decomposing the space into the flows in $\mK$ of the curves of successive intersections with the entry/exit surfaces $E_i$ and $S_i$.  
    
  \begin{prop}\cite[Proposition~10.1]{HR2016a}  \label{prop-levels}
There is a well-defined level function 
\begin{equation}
n_0 \colon \fM_0 \to \mN = \{0,1,2,\ldots\} \ , 
\end{equation}
where the preimage $n_0^{-1}(0) = \tau(\cR')$,  the preimage
$n_0^{-1}(1)$ in the union of two infinite propellers which are
asymptotic to $\tau(\cR')$, and for $\ell  > 1$  the preimage
$n_0^{-1}(\ell)$ is a countable union of finite propellers.
\end{prop}
     The precise description of propellers, both finite and infinite, is given in   \cite[Chapters 11, 12]{HR2016a}, and the decomposition is  made precise there. We give a general sketch of the idea.

\begin{figure}[!htbp]
\centering
\begin{subfigure}[c]{0.4\textwidth}{\includegraphics[height=46mm]{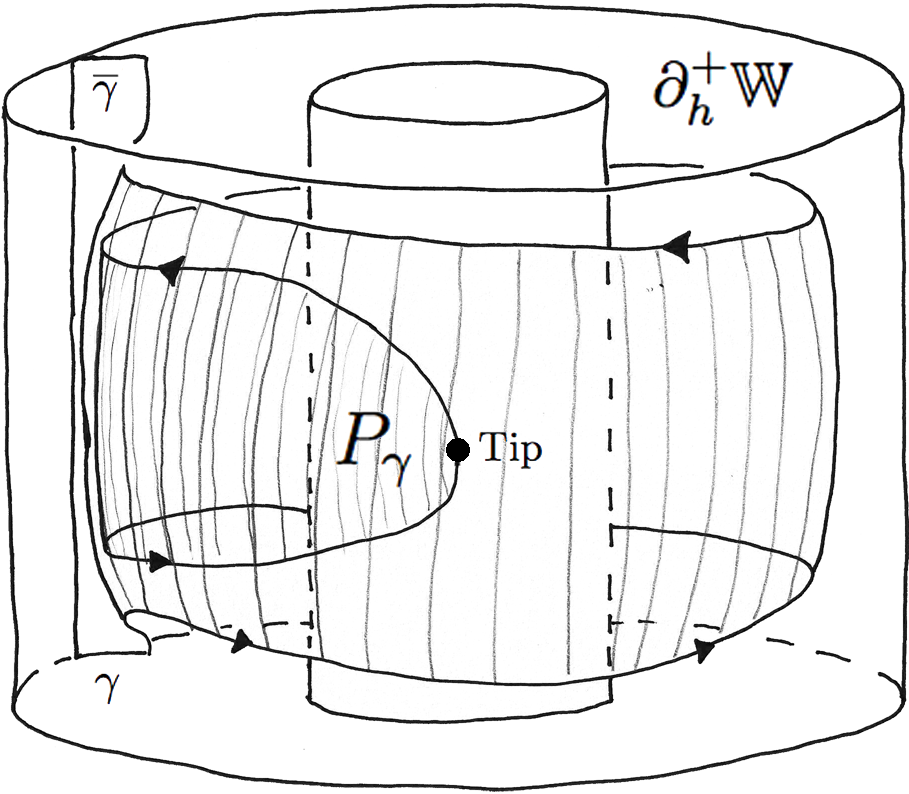}}\end{subfigure}
\begin{subfigure}[c]{0.5\textwidth}{\includegraphics[height=38mm]{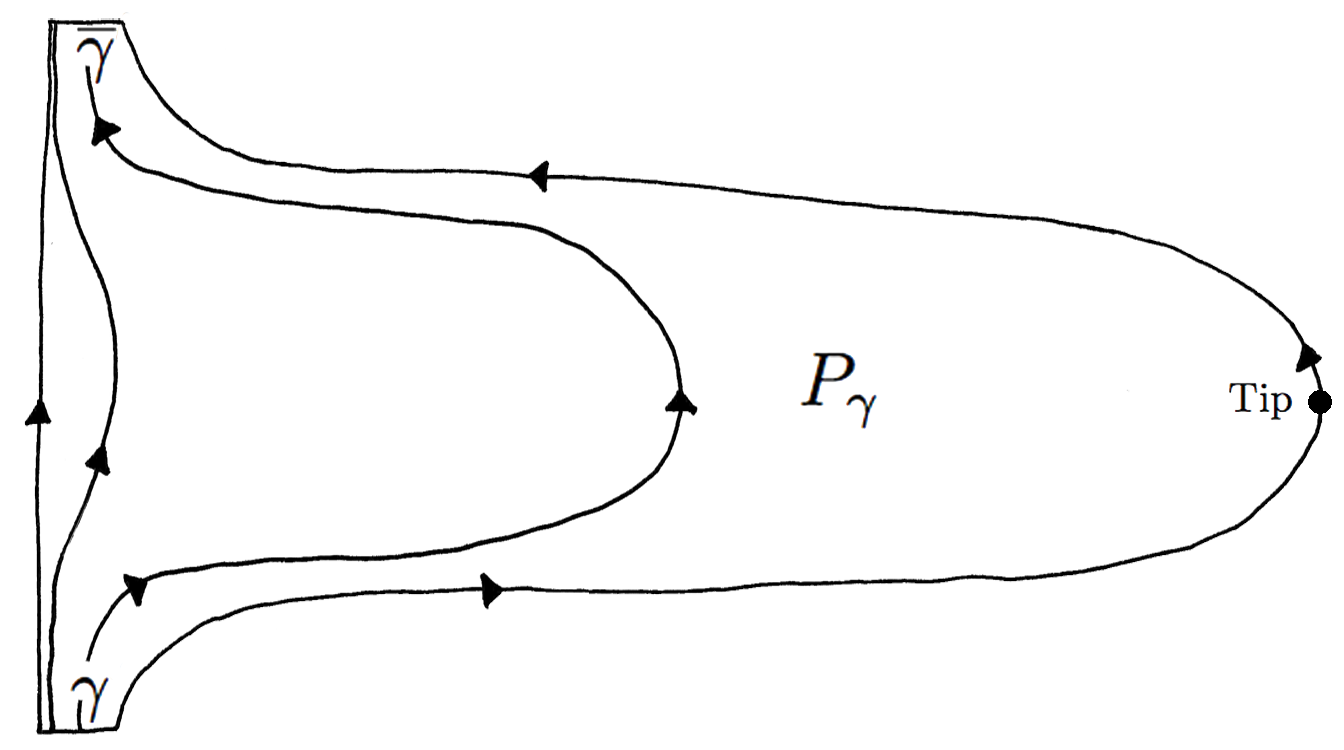}}\end{subfigure}
\caption{\label{fig:propeller}   Embedded and flattened finite propeller}
\vspace{-6pt}
\end{figure}

A {\it propeller} is  an embedded surface in $\mW$  that results from the Wilson flow $\Psi_t$ of a curve $\gamma \subset \partial_h^- \mW$ in the bottom face of $\mW$. Such a surface has   the form of a ``tongue'' wrapping around the core cylinder $\cC$.
  Figure~\ref{fig:propeller}  illustrates a ``typical'' finite propeller $P_{\g}$
  as  a compact  ``flattened'' propeller on the right, and   its
  embedding in $\mW$ on the left. Observe that for any $x\in \g$ the radius of $x$ is strictly bigger
  than 2.  An
  infinite propeller is not closed, and its boundary curve is the
  orbit of an entry point with radius $2$, hence  limits on the Reeb cylinder $\cR$. The embedding of an infinite propeller is highly dependent on the shape of the curve $\gamma$ near the cylinder $\cC$, and on the dynamics of the Wilson flow near its periodic orbits.

\begin{figure}[!htbp]
\centering
{\includegraphics[width=110mm]{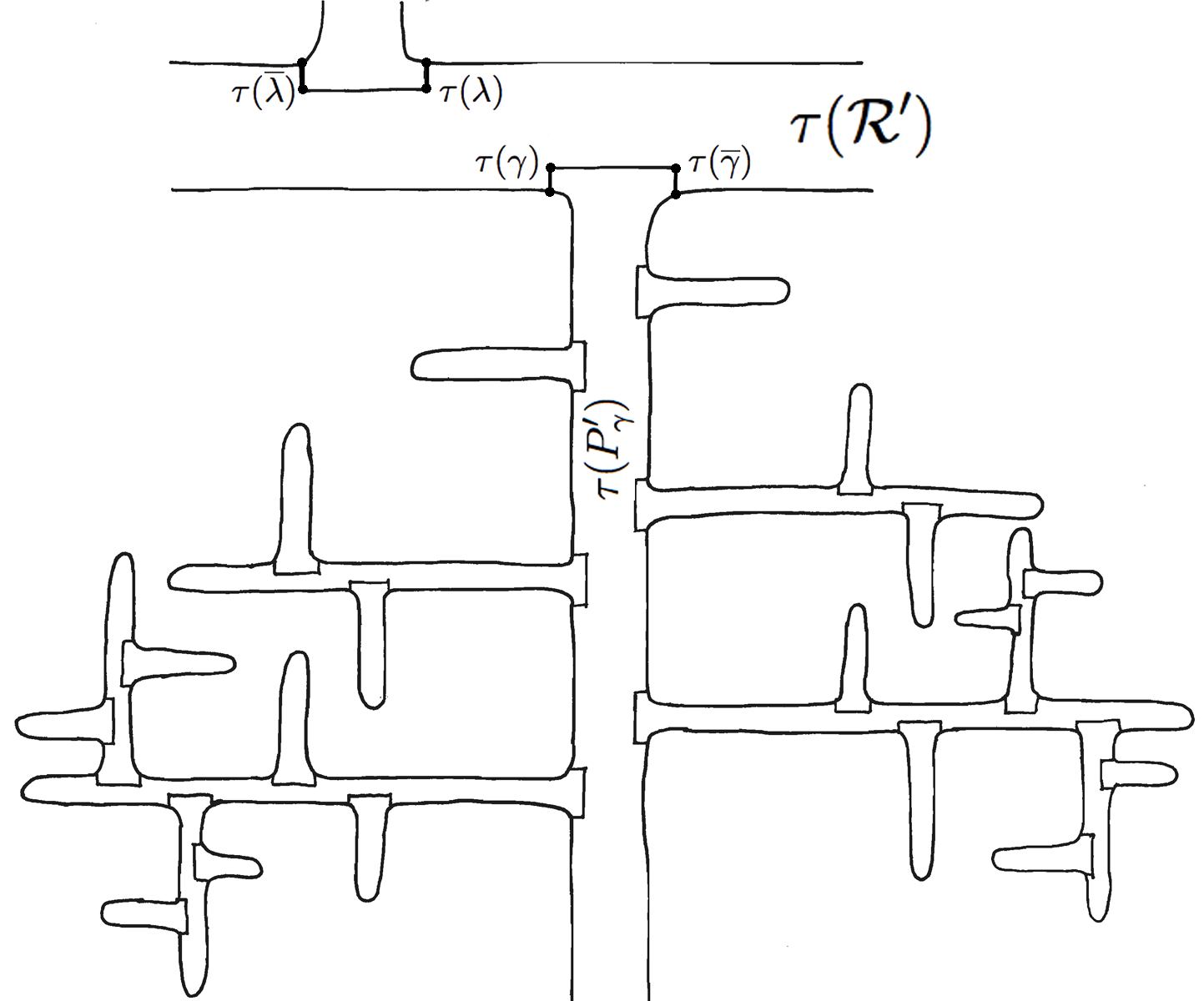}}
\caption{\label{fig:choufleur}  Flattened part of $\fM_0$}
\vspace{-6pt}
\end{figure}

Figure~\ref{fig:choufleur} gives a model for   $\fM_0$, though the
distances along propellers are not to scale, and there is a hidden
simplification in that there may be ``bubbles'' in the surfaces which
are suppressed in the illustration. A bubble is a branching surface attached to the
interior regions of a propeller,   and are analyzed in Chapters~15 and 18 of  \cite{HR2016a}. 
Also, all the propellers represented in Figure~\ref{fig:choufleur} have roughly the same width when embedded in $\mK$, which is the width  of the Reeb cylinder.
   
We make some further comments on the properties of $\fM_0$ as illustrated in Figure~\ref{fig:choufleur}. 
The upper horizontal band the figure represents the notched Reeb cylinder. The flow of the special point $p_1^-$ is the curve along the bottom edge of the image $\tau(\cR')$. When the flow crosses the curve $\tau(\gamma) = \tau(\cR') \cap E_1$, it turns to the right and enters the infinite propeller at level $1$, and follows the left edge of this vertical strip downward, along the Wilson flow of a point with $r=2$ until it intersects the secondary entry surface $E_1$ again. It then turns to the right in the flow direction, and enters a finite propeller at level $2$. In the case pictured, it then flows upward until it crosses the annulus $\cA = \{z=0\} \subset \mW$, which corresponds to the tip of the propeller. It then reverses direction and flows until it crosses the secondary exit surface $S_1$, and resumes flowing downward along the infinite level $1$ propeller. However, as this is following a Wilson orbit, the $z$-values of this part of the orbit are increasing towards $-1$. 

This procedure continues repeatedly, though as the curve moves further
down the level $1$ propeller, the $z$-values get closer to $-1$, and
hence the flow in the side level 2 propellers intersects the secondary
entry region $E_1$ increasingly often, before flowing through the tip
of the corresponding level 2 propeller, and reversing its march through either a secondary exit surface $S_i$ or another secondary entry face $E_i$. This process can be viewed as a geometric model for the recursive description of the flow dynamics  as described using programming language in   \cite[Section~5]{Kuperbergs1996}. A key point is that the lengths of the side branches, while finite, increase in length and branching complexity as the orbit moves downwards along the vertical level $1$ propeller. 
A similar scenario plays out when following the upper infinite propeller, whose initial segment is all that is   illustrated in Figure~\ref{fig:choufleur}.

 The   two infinite propellers which constitute $n_0^{-1}(1)$ are well-understood, but the   finite propellers which constitute the sets  $n_0^{-1}(\ell)$  for $\ell > 1$,  pictured as the side branching surfaces in  Figure~\ref{fig:choufleur},  these may   defy a systematic description without imposing some form of generic hypotheses on the construction of the flow.
    
   On the other hand, for a generic Kuperberg flow as in Definition~\ref{def-generic}, the work   
    \cite{HR2016a}  gives  a reasonably complete description of the components of the level decomposition of $\fM_0$.
 These results are used to show:

   \begin{thm} \cite[Theorem~19.1]{HR2016a} \label{thm-zippered} 
If $\Phi_t$ is  a generic Kuperberg flow on $\mK$, then $\fM$ is a \emph{zippered lamination}. 
 \end{thm}
The definition of a zippered lamination is technical, and given in \cite[Definition~19.3]{HR2016a}. The notion can be summarized by the conditions that $\fM$ is a union of $2$-dimensional submanifolds of $\fM$, and admits a finite cover by special foliation charts which are maps of   subsets of $\fM$ to a measurable product of  a disk with boundary in $\mR^2$ with a Cantor set. In particular, this covering property enables the construction of  the transverse holonomy maps along the leaves of the lamination $\fM$. 
 The structure of the submanifold $\fM_0$ is key to understanding the entropy invariants of the flow, and also conjecturally the Hausdorff dimensions of its closed invariant sets, as will be discussed further in Section~\ref{sec-entropy}.

  \section{Denjoy Theory for laminations}\label{sec-denjoy}

The basic problem in the study of the dynamics of Kuperberg flows is to understand  the topological and ergodic structure of its minimal set $\Sigma$. We have $\Sigma \subset \fM$, and   the most basic problem is the following:

\begin{prob}\label{prob-minset}
Give conditions on a Kuperberg flow which imply that $\Sigma = \fM$.
\end{prob}

The equality $\Sigma = \fM$ is a remarkable conclusion, as the flow of the special orbits $p_i^{\pm} \in \mK$ constitute the boundary of the submanifold $\fM_0$ which is dense in $\fM$, so this asserts that the boundaries of the path connected components of $\fM_0$ are dense in the space itself! This property seems highly improbable. However,    \cite[Th\'eor\`eme, page 302]{Ghys1995}  states that there exists Kuperberg flows for which $\Sigma = \fM$, and hence the minimal set $\Sigma$ is $2$-dimensional.  The result  \cite[Theorem~17]{Kuperbergs1996} gives an explicit analytic    flow for which  $\Sigma = \fM$.

 The idea behind these examples is based on the observation that the
 orbit   $\{\Phi_t(p_1^-) \mid -\infty < t < \infty\}$   of a special
 point $p_1^-$ contains the boundary of all the level 2 propellers
 represented in Figure~\ref{fig:choufleur}, thus it contains the tips
 of these propellers. As the level 2 propellers get longer, the tips
 have smaller radius that tends to 2. The points corresponding to the
 tips are contained in the annulus $\tau(\cA)=\tau(\{z=0\})$, and thus
 accumulate on the Reeb cylinder $\tau(\cR')$.
  The proof of the following
   result was inspired by the proof of
   \cite[Theorem~17]{Kuperbergs1996}, and uses these ideas  to show:

 \begin{thm}\cite[Theorem~17.1]{HR2016a} \label{thm-denjoy}
Let $\Phi_t$ be a generic Kuperberg flow on $\mK$,  then   
$\Sigma  = \fM$.
\end{thm}
The proof of Theorem~\ref{thm-denjoy} uses the generic hypotheses on both the Wilson flow and the insertion maps, to obtain estimates on the density of the orbit 
 $\{\Phi_t(p_1^-) \mid -\infty < t < \infty\}$ near to $\tau(\cR')$. While the calculations in \cite{HR2016a} use the quadratic assumptions on the maps, it seems reasonable to expect that the calculations also work with suitable modifications for the case when 
 the estimates have higher order approximations.
 
\begin{prob}\label{prob-ndensity}
Let $\Phi_t$ be a Kuperberg flow on $\mK$ which satisfies  Hypothesis~\ref{hyp-genericWn} for some even $n \geq 2$,    Hypothesis~\ref{hyp-polynomialn} for some possibly different value of $n \geq 2$, and otherwise satisfies the generic hypotheses. Show that $\Sigma = \fM$.
 \end{prob}
 
 One way to ensure that the hypotheses of Problem~\ref{prob-ndensity} are satisfied 
 is to assume the construction is analytic. 
 
 \begin{prob}\label{prob-adensity}
Let $\Phi_t$ be an analytic Kuperberg flow on $\mK$.  Show that $\Sigma = \fM$.
 \end{prob}

  The other possibility is that these two invariant sets are distinct.   Theorem~19 of \cite{Kuperbergs1996} constructs a piecewise-linear (PL) Kuperberg flow such that the minimal set $\Sigma$ is $1$-dimensional, and thus the inclusion $\Sigma \subset \fM$ is proper. They also assert that there are examples of PL-flows for which the minimal set is $2$-dimensional.  
      
  \begin{prob}\label{prob-undensity}
Let $\Phi_t$ be a Kuperberg flow on $\mK$, smooth or possibly only $C^1$ or Lipschitz. Find conditions on the construction which ensure that  the minimal set $\Sigma$ is $1$-dimensional. For example, if the insertion maps do not satisfy Hypothesis~\ref{hyp-polynomialn}, and in fact are infinitely flat at the special point, is it possible that the inclusion $\Sigma \subset \fM$ is proper? 
 \end{prob}

There is one other aspect of the relationship between $\Sigma$, $\fM$ and the non-wandering set $\Omega$ to mention. The following result  is a direct consequence of Theorem~\ref{thm-minimal} above. 
\begin{thm}\label{thm-wander2}
Let $\Phi_t$ be a Kuperberg flow of $\mK$. Then  $\Sigma \subset \Omega    \subset \fM$.
\end{thm}

The proof of the following result uses  the assumption that the flow is generic to obtain in   \cite[Chapter~16]{HR2016a} a structure theory for the wandering set $\fW$, and hence  to conclude:
 \begin{thm}\cite[Theorem~1.3]{HR2016a} \label{thm-nonwandering}
Let $\Phi_t$ be a generic Kuperberg flow on $\mK$,  then   $\Sigma =   \Omega = \fM$.
\end{thm}
The key point in the proof of this result is to analyze the points in
the complement $\mK - \fM$,  and show that their orbits must include points in the set $\{x \in \fM \mid r(x) < 2\}$, which implies the equality. 

  We mention another natural problem concerning analytic Kuperberg flows, for which one expects additional dynamical properties to be true.
  
    \begin{prob}\label{prob-avsn}
Let $\Phi_t$ be a Kuperberg flow on $\mK$.  Find dynamical properties of $\Phi_t$  which distinguish the cases where the construction is real analytic   from the smooth (possibly non-generic) case.
 \end{prob}

The study of the relationship between   $\Sigma$ and $\fM$ suggests considering a more general question, which is a type of Denjoy theorem for $2$-dimensional laminations, or matchbox manifolds in the terminology of \cite{ClarkHurder2013}.

\begin{prob}\label{prob-denjoy}
Let $\cL$ be a compact connected $2$-dimensional, codimension 1,
lamination, possibly with boundary, and let $\cX$ be a smooth vector
field tangent to the leaves of $\cL$. If the boundary is non-empty, we assume that $\cX$ is tangent
to the boundary.  If $\cL$ is minimal and the flow of  $\cX$ has no periodic orbits, show that every orbit is dense.
 \end{prob}

 The question is whether the equality $\Sigma = \fM$ for Kuperberg flows might follow from a more general ``Denjoy Principle'' which is independent of the embedding of the space $\fM \subset \mK$.  For example, can the proof of the traditional Denjoy Theorem for $C^2$-flows on the $2$-torus $\mT^2$ be adapted to work for laminations? If so, what are the minimal hypotheses required to obtain such a result?

  \section{Growth, slow entropy,  and Hausdorff dimension}\label{sec-entropy}
 
We next consider invariants of Kuperberg flows derived from the choice of  a Riemannian metric on  $\mK$.
These include the area growth rate of the embedded surface $\fM_0 \subset \mK$, the slow entropy of the flow $\Phi_t$ on $\mK$, and the Hausdorff dimensions of the closed invariant sets $\Sigma$ and $\fM$. The work \cite{HR2016a} contains results on these properties for generic flows, but almost nothing is known about them for the case of non-generic flows.

 Choose a Riemannian metric on $\mK$, then 
 the smooth embedded submanifold $\fM_0 \subset \mK$ with boundary inherits a Riemannian metric. Let  $d_{\fM}$ denote the associated path-distance function on $\fM_0$.
 Fix  the basepoint $\omega_0 = (2, \pi, 0) \in \tau(\cR')$ and  let 
 $\ds B_{\omega_0}(s)  = \{ x \in \fM_0 \mid d_{\fM}(\omega_0 , x) \leq s\}$  
be the closed ball of radius $s$ about the basepoint $\omega_0$. 
Let $\A(X)$ denote the Riemannian area of a Borel subset $X \subset \fM_0$. Then $\mathrm{Gr}(\fM_0, s) = \A(B_{\omega_0}(s))$ is called the \emph{growth function} of $\fM_0$.

Given functions $f_1, f_2 \colon [0,\infty) \to [0, \infty)$, we say that $f_1 \lesssim f_2$ if there exists constants $A, B, C > 0$ such that for all $s \geq 0$, we have that $  f_2(s) ~ \leq ~ A \cdot f_1(B \cdot s) + C$.
 Say that    $f_1 \sim f_2$ if both $f_1 \lesssim f_2$ and $f_2 \lesssim f_1$ hold.    This defines   equivalence relation on functions, which defines their   \emph{growth type}.

The growth function $\mathrm{Gr}(\fM_0, s)$ for $\fM_0$ depends upon the choice of Riemannian metric on $\mK$ and basepoint $\omega_0 \in \fM_0$, however  the growth type of $\mathrm{Gr}(\fM_0, s)$  is   independent of   these choices,.

We say that $\fM_0$ has \emph{exponential growth type} if $\mathrm{Gr}(\fM_0, s) \sim \exp(s)$. Note that $\exp(\lambda \, s) \sim \exp(s)$ for any $\lambda > 0$, so there is only one growth class of ``exponential type''.
We say that $\fM_0$ has \emph{nonexponential growth type} if $\mathrm{Gr}(\fM_0, s)  \lesssim \exp(s)$ but $\exp(s) \not\lesssim \mathrm{Gr}(\fM_0, s)$. 
We also have the subclass of   nonexponential growth type,  where $\fM_0$ has  \emph{quasi-polynomial growth type} if there exists $d \geq 0$ such that $\mathrm{Gr}(\fM_0, s)  \lesssim s^d$.  
The growth type of a leaf of a foliation or lamination is an
entropy-type invariant of its dynamics, as discussed in
\cite{Hurder2014}. 

For an embedded propeller $P_{\gamma} \subset \mK$ the area of the
propeller increases as it makes successive revolutions around the core
cylinder, as illustrated in Figure~\ref{fig:propeller}, and this
increase is proportional, with uniform bounds above and below, to the
number of revolutions times the area of the Reeb cylinder $\cR$.
Thus, the growth type of $\mathrm{Gr}(\fM_0, s)$ is a measure of the
number of branches  and their length in $\fM_0$ within a given distance $s$ from $\omega_0$ along the surface. It is thus a measure of the complexity of the recursive procedure which is used in the level decomposition of $\fM_0$. 

\begin{prob}\label{prob-growthtype1}
 Show that the growth type of $\mathrm{Gr}(\fM_0, s)$ for a Kuperberg
 flow is always  nonexponential. 
\end{prob}

This problem was answered in \cite{HR2016a} in the case where the flow
is generic. Under the additional hypothesis on the insertion maps
$\sigma_i$ for $i=1,2$, which is that they have ``slow growth'', the
following result is proved.
 
\begin{thm} \cite[Theorem~22.1]{HR2016a}\label{thm-volumegrowth}
Let $\Phi_t$ be a generic Kuperberg flow. If the insertion maps $\sigma_i$ for $i=1,2$  have ``slow growth'',  then the growth type  of $\fM_0$ is nonexponential, and satisfies
$\ds \exp(\sqrt{s})   \lesssim \mathrm{Gr}(\fM_0, s)$.  In particular, $\fM_0$ does  not have quasi-polynomial growth type.
\end{thm}

The definition of slow growth is given in  \cite[Definition~21.11]{HR2016a}, and will not be recalled here,  as it requires some background preparations. Also defined in that work is the notion of ``fast growth'' in \cite[Definition~21.12]{HR2016a}

The previous theorem suggests two questions:
\begin{prob}\label{prob-growthtype2}
 Show that the growth type of $\mathrm{Gr}(\fM_0, s)$ for a \emph{generic} Kuperberg flow whose insertion maps have slow growth is precisely the growth type of  the function $\ds \exp(\sqrt{s})$.
\end{prob}
It seems reasonable to expect   this problem has a positive answer, or especially in the case where the flow is also analytic. 
The following problem is more open-ended, and   likely much more difficult.
\begin{prob}\label{prob-growthtype3}
How does the growth type of $\mathrm{Gr}(\fM_0, s)$ for a   Kuperberg flow depend on the geometry of the insertion maps, and the germ of the Wilson vector field in a neighborhood of the periodic orbits?
\end{prob}

Part of the motivation for the study of the growth function $\mathrm{Gr}(\fM_0, s)$ is its relation to the topological entropy invariants for the flow $\Phi_t$. 
We define the entropy invariants of the flow $\Phi_t$  using a variation of the Bowen formulation of   topological entropy  \cite{Bowen1971,Walters1982} for a flow $\vp_t$  on a compact metric space $(X, d_X)$, which    is symmetric in the role of the time variable $t$.    For a flow $\vp_t$ on $X$, and for $\e > 0$, two points $p,q\in X$ are said to be 
\emph{$(\vp_t , T, \e)$-separated} if 
\begin{equation}\label{eq-separated}
d_X(\vp_t(p),\vp_t(q))>\e \quad \mbox{for some} \quad -T\leq t\leq T ~ .
\end{equation}
A set $E \subset X$ is \emph{$(\vp_t , T, \e)$-separated} if  all pairs of distinct 
points in $E$ are $(\vp_t , T, \e)$-separated. Let $s(\vp_t , T, \e)$ be the maximal
cardinality of a $(\vp_t , T, \e)$-separated set in $X$. The growth type of the function $s(\vp_t , T, \e)$ is called the $\e$-growth type of $\vp_t$, and we can then study the behavior of the growth type as $\e \to 0$.

The topological entropy of the flow $\vp_t$ is then defined by
\begin{equation}\label{eq-topentropy}
h_{top}(\vp_t)= \frac{1}{2} \cdot \lim_{\e\to  0} \left\{ \limsup_{T\to\infty}\frac{1}{T}\log(s(\vp_t , T, \e)) \right\} \ .
\end{equation}
Moreover, for a compact space $X$, the entropy $h_{top}(\vp_t)$ is independent of the choice of metric $d_X$.
 
A relative form of the topological entropy for a flow $\vp_t$ can be defined for any   subset $Y \subset X$, by requiring that the collection of distinct $(\vp_t , T, \e)$-separated  points used in the definition \eqref{eq-separated}  be contained  in $Y$. The restricted topological entropy $h_{top}(\vp_t | Y)$ is bounded above by $h_{top}(\vp_t)$. 

The notion of slow entropy was introduced in the papers  \cite{deCarvalho1997,KatokThouvenot1997}, and there is the related notion of the entropy dimension \cite{DHP2011}, given as follows:

\begin{defn}\label{def-slowentropy}
For a flow $\vp_t$ on $X$, and $\alpha > 0$, the $\alpha$-slow entropy of $\vp_t$ is given by 
\begin{equation}\label{eq-slowentropy}
h^{\alpha}_{top}(\vp_t)= \frac{1}{2} \cdot \lim_{\e\to  0} \left\{ \limsup_{T\to\infty}\frac{1}{T^{\alpha}}\log(s(\vp_t , T, \e)) \right\} \ .
\end{equation}
\end{defn}

\begin{defn}\label{def-entropydim}
For a flow $\vp_t$ on $X$,   the  entropy dimension of $\vp_t$ is given by 
\begin{equation}\label{eq-entropydim}
 {\rm{Dim}}_{h}(\vp_t)=  \inf_{\alpha > 0} \left\{  h^{\alpha}_{top}(\vp_t) \right\}  = 0\ .
\end{equation}
\end{defn}
For a smooth flow on a compact manifold, we have $0 \leq {\rm{Dim}}_{h}(\vp_t) \leq 1$.

   Katok   proved in \cite[Corollary~4.4]{Katok1980} that for a   $C^2$-flow $\vp_t$ on  a compact $3$-manifold, its topological entropy $h_{top}(\vp_t)$ is bounded above by   the exponent of the  rate of growth of its periodic orbits.  In particular, Katok's result can be applied to an aperiodic  flow obtained by inserting Kuperberg plugs, and  it follows that:
\begin{thm} 
Let $\Phi_t$ be a Kuperberg flow, then the restricted entropy    $h_{top}(\Phi_t | \fM) = 0$.
\end{thm}

Using the choice of a rectangle $\bRt \subset \mK$ to a generic Kuperberg flow, we associate a pseudogroup $\cGM$ formed by the return maps to $\bRt$, acting on the transverse Cantor set $\fC$ to the intersection $\fM \cap \bRt$.  This is described in Chapter~21 of \cite{HR2016a}. There is a notion of $\alpha$-slow entropy $h_{GLW}^{\alpha}(\cGM)$  associated to this pseudogroup, defined by \cite[Formula (165)]{HR2016a} which is a variation on the entropy for pseudogroups introduced in \cite{GLW1988}. Then it was show there that:

\begin{thm}\cite[Theorem~21.10]{HR2016a} \label{thm-psgslowentropy}
Let $\Phi_t$ be a generic Kuperberg flow. If the insertion maps $\sigma_j$  have ``slow growth'', then $h_{GLW}^{1/2}(\cGM) > 0$, and thus the entropy dimension of the pseudogroup action on $\fC$ is bounded below by $1/2$.
\end{thm}

It is natural to ask if the ``slow growth'' hypothesis in Theorem~\ref{thm-psgslowentropy} is necessary:
 
\begin{prob}\label{prob-entdim1}
Let $\Phi_t$ be a generic Kuperberg flow. Show that  $\ds {\rm{Dim}}_{h}(\Phi_t) \geq 1/2$.
\end{prob}
 
 There is a more general variation on this problem, which is possibly more precise as well:
 
 \begin{prob}\label{prob-entdim2}
Let $\Phi_t$ be a   Kuperberg flow, and suppose that the growth type of $\fM_0$ is at least that of the function $n^{\alpha}$, for $0 < \alpha < 1$. Show that  $\ds {\rm{Dim}}_{h}(\Phi_t) \geq \alpha$.
\end{prob}

The discussion and proofs in \cite[Chapters 20 and 21]{HR2016a} contain various arguments which support posing these questions, though the material there does not appear to be sufficient to show these two problems have positive solutions. A key aspect of the estimation of the $\e$-separation function  $s(\Phi_t , T, \e)$ for a Kuperberg flow is the rate of approach of the orbits of the Wilson flow to the periodic orbits. The generic hypothesis for the Wilson flow is used to give estimates on this rate, which is the source of the exponent $\alpha = 1/2$ in Theorem~\ref{thm-psgslowentropy}. In the non-generic case, this rate of approach may be much slower, and so it takes a much longer period of time for orbits to separate. This suggest that the following problem has a positive solution.

\begin{prob}\label{prob-entdim=0}
Let $\Phi_t$ be a   Kuperberg flow,  and suppose that the Wilson flow used in its construction is infinitely flat near its periodic orbits. Show that  $\ds {\rm{Dim}}_{h}(\Phi_t) =0$.
\end{prob}

 At the other extreme from the consideration of infinitely flat Wilson flows, one can   consider the entropy invariants for PL-versions of the Kuperberg construction, as in  \cite[Section~8]{Kuperbergs1996}. Then we allow  the Wilson flow   to
have a discontinuity in its defining vector field $\cW$ along the
periodic orbits, an    we can obtain   the special points $\omega_i$
are   hyperbolic attracting  for the map Wilson flow.  In this case,
the following seems likely to be true:
   \begin{prob}\label{prob-hyperbolicentropy}
Let $\Phi_t$ be a  PL Kuperberg flow,  constructed from a Wilson flow for which the periodic orbits are hyperbolic attracting when restricted to the cylinder $\cC$.  Show that  $\ds h_{top}(\Phi_t | \fM) > 0$.
 \end{prob}

In general, it seems likely that the dynamical and ergodic theory properties of PL-versions of the Kuperberg construction will have a much wider range of possibilities, as was suggested in the work  \cite{Kuperbergs1996}.

The last set of metric invariants for a Kuperberg flow  to consider are its dimension properties. 

\begin{prob}\label{prob-HDSigma}
Show that the Hausdorff dimension of the minimal set $\Sigma$ satisfies ${\rm HD}(\Sigma) > 1$. 
\end{prob}

\begin{prob}\label{prob-HDM1}
Show that the Hausdorff dimension of the invariant set $\fM$ satisfies ${\rm HD}(\fM) > 2$. 
\end{prob}

\begin{prob}\label{prob-HDM2}
Let $\Phi_t$ be a generic Kuperberg flow. Show that ${\rm HD}(\Sigma)  \geq 5/2$. 
\end{prob}

\begin{prob}\label{prob-HDM3}
Is it possible to construct a Kuperberg flow, possibly using a PL-construction, such that  ${\rm HD}(\Sigma)$ can assume any value between $2$ and $3$? 
\end{prob}

  \section{Shape theory for the minimal set}\label{sec-shape}
 
   Shape theory studies the topological properties of a topological space $\fZ$  using a form of \v{C}ech homotopy theory.
   The natural framework for the  study of topological properties of spaces such as the minimal set $\Sigma$ of a Kuperberg flow is using   shape theory.  
For example,   Krystyna Kuperberg posed the question whether $\Sigma $ has stable shape? Stable shape is discussed below, and  is about the nicest property one can expect  for a minimal set that is not a compact submanifold. There are other shape properties of these spaces which can be investigated.  The results that are known about their shape properties  are  all for the generic case. 

We first give a brief introduction to the notions of shape theory, and introduce stable shape and the movable conditions so that we can formulate the known results and some problems.

   The definition of  {shape} for a topological space $\fZ$   was
   introduced by Borsuk \cite{Borsuk1968,Borsuk1975}. Later
   developments and   results of shape theory are  discussed  in the
   texts \cite{DydakSegal1978,MardesicSegal1982} and the historical
   essay  \cite{Mardesic1999}. See also   the works of
   Fox~\cite{Fox1972} and Morita~\cite{Morita1975}. 

Recall that a \emph{continuum} is a compact, connected metrizable space.
For example, the subspaces $\Sigma$ and $\fM$ of $\mK$ are  compact and connected, so are continua.  We discuss below shape theory for continua.
  \begin{defn}\label{def-shapeapprox}
  Let $\fZ \subset X$ be a continuum embedded in a metric space $X$.
  A \emph{shape approximation}   of $\fZ$ is a 
 sequence  $\fU = \{U_{\ell} \mid \ell =1,2,\ldots\}$  satisfying the conditions:
\begin{enumerate}
\item each $U_{\ell}$ is an open neighborhood of $\fZ$ in $X$ which is homotopy equivalent to a compact polyhedron;
\item $U_{\ell +1} \subset U_{\ell}$ for $\ell \geq 1$, and their closures satisfy $\ds \bigcap_{\ell \geq 1} ~ \oU_{\ell} = \fZ$.
\end{enumerate}
 \end{defn} 
 There is a notion of equivalence of shape approximations for continua $\fZ$ and $\fZ'$. In the case where these spaces are embedded in a manifold, the notion of equivalence is discussed  in \cite[Chapter 23]{HR2016a}. Otherwise, any of the sources cited above give  the more general definitions of equivalence of shape approximations. 
\begin{defn}
Let $\fZ \subset X$ be a compact subset of a connected manifold $X$. 
 Then the \emph{shape of $\fZ$ }is defined to be the equivalence class of a shape approximation of $\fZ$ as above. 
\end{defn}   
It is a basic fact of shape theory   that   two homotopy equivalent  continua have the same shape.
Complete details and alternate approaches to defining the shape of a space are   given in  \cite{DydakSegal1978,MardesicSegal1982}.  
An overview of   shape theory for continua embedded in Riemannian manifolds        is   given  in \cite[Section 2]{ClarkHunton2012}.

For the purposes of defining the shape of the spaces $\Sigma$ and $\fM$ for a Kuperberg flow, which are both embedded in $\mK$,   their   shape  can   be defined using    a  shape approximation $\fU$     defined by  a descending chain of open $\e$-neighborhoods in $\mK$ of each set. For example,  the open sets
$\ds U_{\ell} = \{ x \in \mK \mid  d_{\mK}(x, \Sigma) < \e_{\ell}\}$  where we have  $0 <  \e_{\ell +1} < \e_{\ell}$ for all $\ell \geq 1$, and  $\ds \lim_{\ell \to \infty}  \, \e_{\ell}   =   0$, give a shape approximation to $\Sigma$.

Now we define two basic   properties of the shape of a space.

  \begin{defn}\label{def-stableshape}
 A continuum $\fZ$ has \emph{stable shape} if it is \emph{shape equivalent} to   a finite polyhedron. That is, there exists a shape approximation $\fU$    such that each inclusion $\iota \colon U_{\ell +1}\hookrightarrow  U_{\ell}$ induces a homotopy equivalence,  and $U_1$ has the homotopy type of a finite polyhedron. 
 \end{defn}

Some examples of spaces with stable shape are  compact   connected manifolds,  and more generally connected finite $CW$-complexes. A less obvious
 example is the minimal set for a Denjoy flow on $\mT^2$ whose shape is equivalent to the wedge of two circles.  In particular, 
 the minimal set  of an aperiodic $C^1$-flow on plugs as constructed by Schweitzer in \cite{Schweitzer1974} has stable shape. 
 In contrast, the minimal set for a generic Kuperberg flow has very complicated shape, and in particular we have.: 

    \begin{thm}\cite[Theorem~1.5]{HR2016a} \label{thm-unstableshape}
 The  minimal set $\Sigma = \fM$ of a generic Kuperberg flow does  not have stable shape.
 \end{thm}
The proof of this result uses the detailed structure theory for the space $\fM_0$ developed in 
\cite{HR2016a}, to explicitly construct a shape approximation for $\fM$ which is derived from the decomposition of  $\fM_0$ into propellers.  It seems almost certain that with some appropriate additional insights, the following must be true:
\begin{prob}\label{prob-shape2}
Let $\Sigma$ be the minimal set for a Kuperberg flow. Show that $\Sigma$ does not have stable shape.
\end{prob}

There is also an intuitive feeling that the shape of the minimal set depends on the regularity of the flow.
\begin{prob}\label{prob-shape3}
 Let $\Sigma$ be the minimal set for a generic Kuperberg flow. Find shape properties of $\Sigma$ which distinguish it from the minimal set for a non-generic Kuperberg  flow.
\end{prob}

The proof of Theorem~\ref{thm-unstableshape} in \cite[Chapter 23]{HR2016a} uses many of the same properties of the flow which were also used in the calculation that it has non-zero slow entropy. It is natural to speculate   this is not a   coincidence:

\begin{prob}\label{prob-shape4}
Let $\Sigma$ be the minimal set for a Kuperberg flow. Show that the existence of unstable shape approximations to $\Sigma$ implies that the slow entropy $h^{\alpha}_{top}(\Phi_t) > 0$ for some $0 < \alpha <1$.
\end{prob}

A minimal set is said to be exceptional if it is not a submanifold of
the ambient manifold. The previous problem can be stated for any
exceptional minimal set: if the minimal set has unstable shape, must  the slow entropy of the flow positive for some $\alpha$?

There is another, more delicate shape property that can be investigated for the minimal set.

  \begin{defn} \label{def-movable}
 A continuum $\fZ \subset  X$ is said to be \emph{movable  in $X$} if for every neighborhood $U$
 of $\fZ$,  there exists a neighborhood $U_0 \subset U$ of $\fZ$ such
 that, for every neighborhood $W\subset U_0$ of $\fZ$, there is a continuous map $\varphi \colon U_0 \times [0,1] \to U$  satisfying the condition $\varphi(x,0)=x$ and $\varphi(x,1) \in W$ for every point $ x \in U_0$.
\end{defn}

The notion of a movable continuum  was introduced by Borsuk \cite{Borsuk1969} as a generalization of spaces having the shape of an  \emph{absolute neighborhood retract} (ANR's).   See    \cite{ClarkHunton2012,DydakSegal1978,Krasinkiewicz1981,MardesicSegal1982} for further discussions concerning movability.  
 It  is a   subtle problem to construct continuum  which are invariant sets for dynamical systems and which are movable,  but do not have stable shape, such as given in   \cite{Sindelarova2007}.  
Showing the movable property for a space  requires the construction of a homotopy retract  $\varphi$ with the properties stated in the definition, whose existence  can be difficult to achieve in practice. There is an alternate condition on homology groups, weaker than  the movable condition.

 \begin{prop}\label{prop-MLmove}
Let $\fZ$ be a movable continuum with shape approximation $\fU$. Then the    homology groups satisfy   the \emph{Mittag-Leffler Condition}:
 For all $\ell \geq 1$, there exists $p \geq \ell$ such that for any $q \geq p$, the maps on   homology groups for $m \geq 1$ induced by the inclusion maps  satisfy
 \begin{equation}\label{eq-shapehomology}
\text{Image}\left\{ H_m(U_p; \mZ) \to H_m(U_{\ell} ; \mZ) \right\}  = \text{Image}\left\{ H_m(U_q; \mZ) \to H_m(U_{\ell} ; \mZ) \right\}  ~ .
\end{equation}
\end{prop}

 This result is a   special case of a more general   Mittag-Leffler condition, as discussed in detail in \cite{ClarkHunton2012}.  For example,   the above   form of the Mittag-Leffler condition can be used to show that the Vietoris solenoid formed from the inverse limit of coverings of the circle is not movable.
 
We can now state   an additional shape property for the minimal set of a generic Kuperberg flow. 
\begin{thm}\cite[Theorem 1.6]{HR2016a}\label{thm-MLhomology}
Let $\Sigma$ be the minimal set for a generic Kuperberg flow. Then the Mittag-Leffler condition   for homology groups is satisfied. That is, given a shape approximation $\fU = \{U_{\ell}\}$ for $\Sigma$,
then  for any $\ell\geq 1$ there exists $p>\ell$ such that for any $q\geq p$ 
\begin{equation}\label{eq-MLhomology}
Image\{H_1(U_p;\mZ)\to H_1(U_\ell;\mZ)\} = Image \{H_1(U_q;\mZ)\to H_1(U_\ell;\mZ)\}.
\end{equation}
\end{thm}

The proof of Theorem~\ref{thm-MLhomology} in \cite[Chapter 23]{HR2016a} is even more subtle than the proof of Theorem~\ref{thm-unstableshape}, but it suggests the following should be true:

\begin{prob}\label{prob-genericmove1}
Show that   the minimal set $\Sigma$ for a generic Kuperberg flow  is movable.
\end{prob}

On the other hand, it would be very remarkable if the minimal set for all Kuperberg flows should be movable.

\begin{prob}\label{prob-genericmove2}
Construct an example of a Kuperberg flow such that the minimal set $\Sigma$ is not movable.
\end{prob}

  \section{Derived from Kuperberg flows}\label{sec-DK}
  
 The discussions in the previous sections show that there are many choices in the construction of Kuperberg flows, and while all result in aperiodic flows, it is conjectured that many of the other dynamical properties of these flows depend upon the choices made. In this final section, we discuss one further variant on the construction, where the resulting flows are no longer aperiodic. As these constructions use the same method as outlined in Sections~\ref{sec-wilson} and \ref{sec-kuperberg}, this new class of flows are called ``Derived from Kuperberg'' flows, or DK--flows.
  The   DK--flows  were introduced in \cite{HR2016b}, and are constructed by varying the construction of the usual Kuperberg flows so that the periodic orbits are not ``broken open''. Thus, the DK--flows are quite useless as counterexamples to the Seifert Conjecture, but they are obtained by smooth variations of the standard Kuperberg flows, so have interest from the point of view of the properties of Kuperberg flows in the space of flows \cite{CY2015,GY2016,Palis2000,Palis2008,PujalsSambarino2000}.   
The work \cite{HR2016b} gave constructions of DK--flows which in fact have countably many independent horseshoe subsystems, and thus have positive topological entropy.  In this section, we discuss some of the questions that arise with   the construction of DK--flows and the study of their properties.

The DK--flows were constructed at the end of
Section~\ref{sec-kuperberg}, and the generic hypotheses on these flows
was formulated in Definition~\ref{def-genericDK}.  Then in
\cite[Section~9.2]{HR2016b} the \emph{admissibility condition} was
formulated for these flows. Briefly, this condition is that there
exist a constant $C > 0$, which depends only on the generic Wilson flow used in the construction, so that if the vertical offset in the $z$-coordinate of the vertex of an insertion map $\sigma_i^\e$ is $\delta > 0$, then we assume that the horizontal offset  $\e$ satisfies $0 < \e < C \cdot \sqrt{\delta}$. 
Then we have:
\begin{thm}\cite[Theorem 9.5]{HR2016b}\label{thm-genDKpos} For $\e > 0$, let
  $\Phi_t^\e$ be a generic DK--flow on $\mK_\e$ which satisfies the admissibility condition. Then  $\Phi_t^\e$ has an invariant horseshoe dynamical system, and thus $h_{top}(\Phi_t^\e) > 0$.
\end{thm}
The proof of this result requires the introduction of the pseudogroup on a transversal $\bRt$ which was constructed in 
\cite[Chapter 9]{HR2016a} or \cite[Section 7]{HR2016b}.  The second
reference gives just the bare minimum of details required to prove
Theorem~\ref{thm-genDKpos}, while the first reference has a
comprehensive discussion of the   pseudogroups associated to
Kuperberg flows. These details are not required to formulate the
following problems, though we believe that they are the starting point for seeking their solutions.

The following problem may just require a technical extension of the ideas used in the proof of Theorem~\ref{thm-genDKpos}, though it is also possible that there are novel dynamical problems which arise in the study of it.
\begin{prob}\label{prob-DKpos}
Show that the topological entropy $h_{top}(\Phi_t^\e) > 0$ for a DK--flow $\Phi_t^\e$   on $\mK_\e$ with $\e > 0$.
\end{prob}

There is a variant of this problem which is discussed further in \cite{HR2016d}.
 \begin{prob}\label{conj-genDKpos}
 Let $\Phi_t^\e$ be a generic DK--flow on $\mK_\e$ with $\e>0$, and let
 $\phi_t$ be the flow obtained from the construction of $\Phi_t^\e$ by
 taking   a sufficiently small smooth perturbation of  the function
 $g$ used in the construction of the Wilson flow which removes its
 vanishing points.  Does the flow $\phi_t$ has an invariant horseshoe dynamical system, and thus $h_{top}(\phi_t) > 0$? 
 \end{prob}
 
 In any case, by Theorem~\ref{thm-genDKpos}, and possibly by an
 affirmative answer to Problem~\ref{conj-genDKpos}, there exists smooth families of DK--flows with positive entropy which limit on a given generic Kuperberg flow, which has entropy zero.  The horseshoe dynamics of these flows are shown to exist using the shape approximations introduced in \cite[Chapter 23]{HR2016a} and discussed in Section~\ref{sec-shape} above.  These shape approximations are based on the structure theory for the submanifold $\fM_0$ discussed in Section~\ref{sec-lamination}, and thus are reasonably well understood for the case of generic flows. On the other hand, the dynamical properties   of the horseshoes for the perturbed flow $\Phi_t^\e$ are still unexplored. 
 
\begin{prob}\label{prob-horseshoes1}
Let $\{\Phi_t^\e \mid 0 < \e \leq \e_0\}$ be a family of generic DK--flows on $\mK_\e$ which converge to a generic Kuperberg flow $\Phi_t$ in the $C^{\infty}$-topology of flows.  Study the limiting behavior of the periodic orbits for the invariant horseshoes for the flows $\Phi_t^\e$   as $\e \to 0$? The limits of such periodic orbits converge to a closed current supported on the minimal set $\Sigma$. Describe the currents on $\Sigma$ that arise in this way.
\end{prob}

 The construction of the horseshoe dynamics for the generic DK--flows in Theorem~\ref{thm-genDKpos} are based on choosing appropriate compact branches of the embedded surface $\fM_0$ as discussed in \cite{HR2016b}.  These compact surfaces   approximately generate the orbits which define the horseshoe, and so behave much like a template for the horseshoes created \cite{Williams1974,Williams1979}. 
 \begin{prob}\label{prob-templates1}
Show that the horseshoes for a generic DK--flow with positive entropy are carried by templates derived from the compact pieces of $\fM_0$.
\end{prob}
 
  The dynamics of the DK flows appears to be  reminiscent of the analysis of the dynamics of Lorenz attractors, as discussed for example in the survey by Ghys \cite{Ghys2013}. 
Moreover, the variation of the horseshoes for a smooth family for generic DK--flows with positive entropy suggests a comparison with the degeneration in the dynamics of the Lorenz attractors as studied by  de Carvalho and Hall 
\cite{deCarvalho1999,deCH2002a,deCH2002b,Hall1994}.   The analogy between the dynamics of a family of generic DK--flows and a family of Lorenz attractors suggests that the topic is worth further investigation.

  \vfill
  \eject

\end{document}